
%
%
%

\magnification=\magstep1
\baselineskip =5mm
\lineskiplimit =1.0mm
\lineskip =1.0mm

\long\def\comment#1{}

\long\def\blankout #1\eb{}
\def\noblankout{\def\blankout{}\def\eb{}}


\let\properlbrack=\lbrack
\let\properrbrack=\rbrack
\def\ordcomma{,}
\def\ordcolon{:}
\def\ordsemicolon{;}
\def\ordleftparen{(}
\def\ordrightparen{)}
\def\ordleftbrack{\properlbrack}
\def\ordrightbrack{\properrbrack}
\def\rmcomma{\ifmmode ,\else \/{\rm ,}\fi}
\def\rmcolon{\ifmmode :\else \/{\rm :}\fi}
\def\rmsemicolon{\ifmmode ;\else \/{\rm ;}\fi}
\def\rmleftparen{\ifmmode (\else \/{\rm (}\fi}
\def\rmrightparen{\ifmmode )\else \/{\rm )}\fi}
\def\rmleftbrack{\ifmmode \properlbrack\else \/{\rm \properlbrack}\fi}
\def\rmrightbrack{\ifmmode \properrbrack\else \/{\rm \properrbrack}\fi}
\catcode`,=\active 
\catcode`:=\active 
\catcode`;=\active 
\catcode`(=\active 
\catcode`)=\active 
\catcode`[=\active 
\catcode`]=\active 
\let,=\ordcomma
\let:=\ordcolon
\let;=\ordsemicolon
\let(=\ordleftparen
\let)=\ordrightparen
\let[=\ordleftbrack
\let]=\ordrightbrack
\let\lbrack=\ordleftbrack
\let\rbrack=\ordrightbrack
\def\rmpunctuation{
\let,=\rmcomma
\let:=\rmcolon
\let;=\rmsemicolon
\let(=\rmleftparen
\let)=\rmrightparen
\let[=\rmleftbrack
\let]=\rmrightbrack
\let\lbrack=\rmleftbrack
\let\rbrack=\rmrightbrack}

\def\writemonth#1{\ifcase#1
\or January\or February\or March\or April\or May\or June\or July%
\or August\or September\or October\or November\or December\fi}

\newcount\mins
\newcount\minmodhour
\newcount\hour
\newcount\hourinmin
\newcount\ampm
\newcount\ampminhour
\newcount\hourmodampm
\def\writetime#1{%
\mins=#1%
\hour=\mins \divide\hour by 60
\hourinmin=\hour \multiply\hourinmin by -60
\minmodhour=\mins \advance\minmodhour by \hourinmin
\ampm=\hour \divide\ampm by 12
\ampminhour=\ampm \multiply\ampminhour by -12
\hourmodampm=\hour \advance\hourmodampm by \ampminhour
\ifnum\hourmodampm=0 12\else \number\hourmodampm\fi
:\ifnum\minmodhour<10 0\number\minmodhour\else \number\minmodhour\fi
\ifodd\ampm p.m.\else a.m.\fi
}

\font\tenrm=cmr10
\font\smallcaps=cmcsc10
\font\eightrm=cmr8
\font\ninerm=cmr9
\font\sixrm=cmr6
\font\eightbf=cmbx8
\font\sixbf=cmbx6
\font\eightit=cmti8
\font\eightsl=cmsl8
\font\eighti=cmmi8
\font\eightsy=cmsy8
\font\eightex=cmex10 at 8pt
\font\sixi=cmmi6
\font\sixsy=cmsy6
\font\ninesy=cmsy9
\font\seventeenrm=cmr17
\font\bigseventeenrm=cmr17 scaled \magstep 1
\font\twelverm=cmr10 scaled \magstep2
\font\seventeeni=cmmi10 scaled \magstep3
\font\twelvei=cmmi10 scaled \magstep2
\font\seventeensy=cmsy10 scaled \magstep3
\font\twelvesy=cmsy10 at 12pt
\font\seventeenex=cmex10 scaled \magstep3
\font\seventeenbf=cmbx10 scaled \magstep3
\catcode`@=11
\def\eightbig#1{{\hbox{$\textfont0=\ninerm\textfont2=\ninesy
\left#1\vbox to6.5pt{}\right.\n@space$}}}
\catcode`@=12
\def\eightpoint{\eightrm \normalbaselineskip=4.5 mm%
\textfont0=\eightrm \scriptfont0=\sixrm \scriptscriptfont0=\fiverm%
\def\rm{\fam0 \eightrm}%
\textfont1=\eighti \scriptfont1=\sixi \scriptscriptfont1=\fivei%
\def\mit{\fam1 } \def\oldstyle{\fam1 \eighti}%
\textfont2=\eightsy \scriptfont2=\sixsy \scriptscriptfont2=\fivesy%
\def\cal{\fam2 }%
\textfont3=\eightex \scriptfont3=\eightex \scriptscriptfont3=\eightex%
\def\bf{\fam\bffam\eightbf} \textfont\bffam\eightbf
\scriptfont\bffam=\sixbf \scriptscriptfont\bffam=\fivebf
\def\it{\fam\itfam\eightit} \textfont\itfam\eightit
\def\sl{\fam\slfam\eightsl} \textfont\slfam\eightsl
\let\big=\eightbig \normalbaselines\rm
\def\caps##1{\nottencaps{##1}}
}
\def\seventeenpoint{\seventeenrm \baselineskip=5.5mm%
\textfont0=\seventeenrm \scriptfont0=\twelverm \scriptscriptfont0=\sevenrm%
\def\rm{\fam0 \seventeenrm}%
\textfont1=\seventeeni \scriptfont1=\twelvei \scriptscriptfont1=\seveni%
\def\mit{\fam1 } \def\oldstyle{\fam1 \seventeeni}%
\textfont2=\seventeensy \scriptfont2=\twelvesy \scriptscriptfont2=\sevensy%
\def\cal{\fam2 }%
\textfont3=\seventeenex \scriptfont3=\seventeenex%
\scriptscriptfont3=\seventeenex%
\def\bf{\fam\bffam\seventeenbf} \textfont\bffam\seventeenbf
}

\def\setheadline #1\\ #2 \par{\headline={\ifnum\pageno=1 
\hfil
\else \eightpoint \noindent
\ifodd\pageno \hfil \caps{#2}\hfil \else
\hfil \caps{#1}\hfil \fi\fi}}

\def\beginsection{} 
\def\datedversion{\footline={\ifnum\pageno=1 \fiverm \hfil
Typeset using plain-\TeX\ on
\writemonth\month\ \number\day, \number\year\ at \writetime{\time}\hfil 
\else \tenrm \hfil \folio \hfil \fi}
\def\tempsetheadline##1{\headline={\ifnum\pageno=1 
\hfil
\else \eightpoint \noindent
\writemonth\month\ \number\day, \number\year,
\writetime{\time}\hfil ##1\fi}}
\def\firstbeginsection##1\par{\bigskip\vskip\parskip
\message{##1}\centerline{\caps{##1}}\nobreak\smallskip\noindent
\tempsetheadline{##1}} \def\beginsection##1\par{\vskip0pt
plus.3\vsize\penalty-250 \vskip0pt plus-.3\vsize\bigskip\vskip\parskip
\message{##1}\centerline{\caps{##1}}\nobreak\smallskip\noindent
\tempsetheadline{##1}}}

\def\finalversion{\footline={\ifnum\pageno=1 \eightrm \hfil 
This paper is in final form.\hfil 
\else \tenrm \hfil \folio \hfil \fi}}

\def\preliminaryversion{\footline={\ifnum\pageno=1 \eightrm \hfil 
Preliminary Version.\hfil 
\else \tenrm \hfil \folio \hfil \fi}}

\def\moreproclaim{\par}
\def\Head #1: {\medskip\noindent{\it #1}:\enspace}
\def\Proof: {\Head Proof: }
\def\Proofof #1: {\Head Proof of #1: }
\def\endproof{\nobreak\hfill$\sqr$\bigskip\goodbreak}
\def\itemi{\item{i)}}
\def\itemii{\item{ii)}}
\def\itemiii{\item{iii)}}
\def\itemiv{\item{iv)}}
\def\itemv{\item{v)}}
\def\itemvi{\item{vi)}}
\def\itemvii{\item{vii)}}
\def\itemviii{\item{viii)}}

\def\dt{\it}
\def\ds#1{{\displaystyle{#1}}}
\def\ts#1{{\textstyle{#1}}}

\def\Abstract\par#1\par{\centerline{\vtop{
\eightpoint
\abovedisplayskip=6pt plus 3pt minus 3pt
\belowdisplayskip=6pt plus 3pt minus 3pt
\moreabstract\parindent=0 true in%
\caps{Abstract}: \ \ #1}}
\abovedisplayskip=12pt plus 3pt minus 9pt
\belowdisplayskip=12pt plus 3pt minus 9pt
\vskip 0.4 true in}
\def\moreabstract{%
\par \hsize = 5 true in \hangindent=0 true in \parindent=0.5 true in}

\def\caps#1{\smallcaps #1}
\def\nottencaps#1{\uppercase{#1}}

\def\firstbeginsection#1\par{\bigskip\vskip\parskip
\message{#1}\centerline{\caps{#1}}\nobreak\smallskip\noindent}

\def\beginsection#1\par{\vskip0pt plus.3\vsize\penalty-250
\vskip0pt plus-.3\vsize\bigskip\vskip\parskip
\message{#1}\centerline{\caps{#1}}\nobreak\smallskip\noindent}

\def\proclaim#1. #2\par{
\medbreak
\noindent{\caps{#1}.\enspace}{\it\rmpunctuation#2\par}
\ifdim\lastskip<\medskipamount \removelastskip
\penalty55\medskip\fi}

\def\Definition: #1\par{
\Head Definition: #1\par
\ifdim\lastskip<\medskipamount \removelastskip
\penalty55\medskip\fi}

\def\Problem #1: #2\par{
\Head Problem #1: #2\par
\ifdim\lastskip<\medskipamount \removelastskip
\penalty55\medskip\fi}


\def\sqr{\vcenter {\hrule height.3mm
\hbox {\vrule width.3mm height 2mm \kern2mm
\vrule width.3mm } \hrule height.3mm }}

\def\references#1{{
\frenchspacing
\eightpoint
\rmpunctuation
\halign{\bf##\hfil & \quad\vtop{\hsize=5.5 true
in\parindent=0pt\hangindent=3mm \strut\rm##\strut\smallskip}\cr#1}}}

\def\ref[#1]{{\bf [#1]}}

\catcode`@=11 
\def\vfootnote#1{\insert\footins\bgroup
\eightpoint
\interlinepenalty=\interfootnotelinepenalty
\splittopskip=\ht\strutbox
\splitmaxdepth=\dp\strutbox \floatingpenalty=20000
\leftskip=0pt \rightskip=0pt \spaceskip=0pt \xspaceskip=0pt
\textindent{#1}\footstrut\futurelet\next\fo@t}

\def\footremark{\insert\footins\bgroup
\eightpoint\it\rmpunctuation
\interlinepenalty=\interfootnotelinepenalty
\splittopskip=\ht\strutbox
\splitmaxdepth=\dp\strutbox \floatingpenalty=20000
\leftskip=0pt \rightskip=0pt \spaceskip=0pt \xspaceskip=0pt
\noindent\footstrut\futurelet\next\fo@t}
\catcode`@=12


\def\Bbb{\bf}
\def\E{{\Bbb E}}
\def\R{{\Bbb R}}
\def\Z{{\Bbb Z}}
\def\N{{\Bbb N}}
\def\C{{\Bbb C}}
\font\specialeightrm=cmr10 at 8pt
\def\R{\hbox{\rm I\kern-2pt R}}
\def\Z{\hbox{\rm Z\kern-3pt Z}}
\def\N{\hbox{\rm I\kern-2pt I\kern-3.1pt N}}
\def\C{\hbox{\rm \kern0.7pt\raise0.8pt\hbox{\specialeightrm I}\kern-4.2pt C}} 
\def\E{\hbox{\rm I\kern-2pt E}}

\def\F{{\cal F}}

\def\invc{{c^{-1}}}
\def\half{{1\over 2}}

\def\qop{{q\over p}}

\def\list#1,#2{#1_1$, $#1_2,\ldots,$\ $#1_{#2}}
\def\lists#1{#1_1$, $#1_2,\ldots}

\def\lnorm{\left\|}
\def\rnorm{\right\|}
\def\normo#1{\lnorm #1 \rnorm}
\def\widedot{\,\cdot\,}
\def\normdot{\normo{\widedot}}

\def\lmod{\left|}
\def\rmod{\right|}
\def\modo#1{\lmod #1 \rmod}

\def\implies{$\Rightarrow$}
\def\iff{$\Leftrightarrow$}


\def\Deltacond{$\Delta_2$-condition}
\def\phifunction{$\varphi$-function}
\def\Nfunction{{\it N}-function}
\def\conditionJ{condition~$(J)$}
\def\conditionL{condition~$(L)$}
\def\em{\mathop{\rm em}\nolimits}
\def\lm{\mathop{\rm lm}\nolimits}
\def\T{{\cal T}}


\noblankout
\finalversion

\setheadline Comparison of Orlicz--Lorentz Spaces\\
             Montgomery-Smith

{
\seventeenpoint
\centerline{\bigseventeenrm Comparison of }
\centerline{\bigseventeenrm Orlicz--Lorentz Spaces}
}
\bigskip\bigskip\medskip
\centerline{\caps{S.J.~Montgomery-Smith}%
\footnote{*}%
{Research supported in part by N.S.F.\ Grant DMS 9001796.}%
}
\smallskip
{
\eightpoint
\centerline{\it Department of Mathematics, University of Missouri,}
\centerline{\it Columbia, MO 65211.}
}
\bigskip\bigskip

{
\eightpoint
\centerline{\it I dedicate this paper to my Mother and Father, who as well as
introducing me to mathematics}
\centerline{\it at an early age, shared the arduous task
of bringing me up.}
}
\bigskip\bigskip

\Abstract

Orlicz--Lorentz spaces provide a common generalization of Orlicz spaces and
Lorentz spaces. They have been studied by many authors, including
Masty\l o, Maligranda, and Kami\'nska. In this paper, we consider the
problem of comparing the Orlicz--Lorentz norms, and establish necessary and
sufficient conditions for them to be equivalent. As a corollary, we give
necessary and sufficient conditions for a Lorentz--Sharpley space to be
equivalent to an Orlicz space, extending results of Lorentz and Raynaud. We
also give an example of a rearrangement invariant space that is not an
Orlicz--Lorentz space.

\footremark{A.M.S.\ (1980) subject classification: 46E30.}

\firstbeginsection 1.\enspace Introduction

The most well known examples of Banach spaces are the $L_p$\ spaces. Their
definition is very well known: if $(\Omega,\F,\mu)$\ is a measure space, and
$1\le p \le \infty$, then for any measurable function $f:\Omega\to\C$, the
$L_p$-norm is defined to be
$$ \normo f_p = 
   \left( \int_\Omega \modo{f(\omega)}^p \,d\mu(\omega) \right)^{1/p} $$
for $p<\infty$, and
$$ \normo f_\infty = \mathop{\hbox{\rm ess sup}}\limits_{\omega\in\Omega}
   \modo{f(\omega)} $$
for $p=\infty$.
Then we define the Banach space $L_p(\Omega,\F,\mu)$\ to be the vector space of
all measurable functions $f:\Omega\to\C$\ for which $\normo f_p$\ is finite.

It is natural to search for generalizations of these $L_p$\ spaces. The first
examples are the {\dt Orlicz spaces}. These were first studied by Orlicz
\ref[O] and Luxemburg \ref[L]. We say that $F:[0,\infty)\to[0,\infty)$\ is an
{\dt Orlicz function\/} if $F$\ is non-decreasing and convex with $F(0)=0$. Now
we define the {\dt Luxemburg norm\/} by 
$$ \normo{f}_F = \inf\left\{\, c :
   \int_\Omega F\bigl(\modo{f(\omega)}/c\bigr) \,d\mu(\omega) \le 1 \,\right\} 
   ,$$
whenever $f$\ is a measurable function, and define the {\dt Orlicz space\/}
$L_F(\Omega,\F,\mu)$\ to be those measurable functions $f$\ for which $\normo
f_F$\ is finite. The Orlicz space $L_F$\ is a true generalization of $L_p$, at
least for $p<\infty$: if $F(t) = t^p$, then $L_F = L_p$\ with equality of norms.

The other examples are the {\dt Lorentz spaces}. These were introduced by
Lorentz \ref[Lo1], \ref[Lo2]. If $f$\ is a measurable function, we define the
{\dt non-increasing rearrangement\/} of $f$\ to be 
$$ f^*(x) = \sup\bigl\{\, t : \mu(\modo f \ge t) \ge x \,\bigr\} .$$ 
If $1\le q <
\infty$, and if $w:(0,\infty)\to(0,\infty)$\ is a non-increasing function, we
define the {\dt Lorentz norm\/} of a measurable function $f$\ to be 
$$ \normo f_{w,q} = \left(\int_0^\infty w(x) f^*(x)^q \,dx\right)^{1/q} .$$ 
We
define the {\dt Lorentz space\/} $\Lambda_{w,q}(\Omega,\F,\mu)$\ to be the space
of those measurable functions $f$\ for which $\normo f_{w,q}$\ is finite. These
spaces also represent a generalization of the $L_p$\ spaces: if $w(x) = 1$\ for
all $0\le x <\infty$, then $\Lambda_{w,p} = L_p$\ with equality of norms.

There is one, rather peculiar, choice of the function $w$\ which turns out to
be rather useful. If $1\le q \le p <\infty$, we define the spaces $L_{p,q}$\ to
be $\Lambda_{w,q}$\ with $w(x) = \qop x^{q/p-1}$. A good reference for a
description of these spaces is Hunt \ref[H]. By a
suitable change of variables, the $L_{p,q}$\ norm may also be defined in the
following fashion: 
$$ \normo f_{p,q} = 
   \left(\int_0^\infty \modo{f^*(x^{p/q})}^q \,dx\right)^{1/q} .$$
Thus $L_{p,p}
= L_p$\ with equality of norms. The reason for this definition is
that for any measurable set $A\in\F$, we have that $\normo{\chi_A}_{p,q} =
\normo{\chi_A}_p = \mu(A)^{1/p} $. Thus $L_{p,q}$\ is a space identical to
$L_p$\ for characteristic functions, but `glued' together in a $L_q$\ fashion.

In all the spaces defined above, if we only desire to study quasi-Banach spaces
rather than Banach spaces, we may remove some of the restrictions placed upon
the defining parameters. Thus with the $L_p$\ spaces, we need only have
$p>0$. With the Orlicz spaces $L_F$\ and the Lorentz space $\Lambda_{w,q}$, we
may weaken the restrictions that $F$\ be convex and that $w$\ be
non-increasing (we omit details). The spaces $\Lambda_{w,q}$\ so obtained were
studied by Sharpley \ref[S], and so we might call then Lorentz--Sharpley spaces.
With the $L_{p,q}$\ spaces, we need only have $0<p<\infty$\ and $0<q\le\infty$,
where if $q=\infty$, we define the Lorentz norm by 
$$ \normo f_{p,\infty} = \sup_{x\ge0} x^{1/p} f^*(x) .$$

Now we come to the object of the paper, the {\dt Orlicz--Lorentz spaces}. These
are a common generalization of the Orlicz spaces and the Lorentz spaces. They
have been studied by Masty\l o (see part~4 of \ref[My]), Maligranda \ref[Ma], and
Kami\'nska \ref[Ka1], \ref[Ka2], \ref[Ka3]. If $G$\
is an Orlicz function, and if $w:[0,\infty)\to[0,\infty)$\ is a non-increasing
function, we define the {\dt Orlicz--Lorentz norm\/} of a measurable function
$f$\ to be 
$$ \normo f_{w,G} = \inf\left\{\, c :
   \int_0^\infty w(x) G\bigl(f^*(x)/c\bigr) \,dx \le 1 \,\right\} 
   .$$
We define the {\dt Orlicz--Lorentz space\/} $\Lambda_{w,G}(\Omega,\F,\mu)$\ to be
the vector space of measurable functions $f$\ for which $\normo f_{w,G}$\ is
finite. If we do not require that the space be a Banach space, but only a
quasi-Banach space, we may weaken the restrictions placed upon $G$\ and $w$\ as
we did for $L_F$\ and $\Lambda_{w,p}$\ above.

We shall not work with this definition of the Orlicz--Lorentz space, however, but
with a different, equivalent definition that bears more resemblance to the spaces
$L_{p,q}$. This definition we give in the following section.

\beginsection 2.\enspace Definitions

First we define \phifunction s. These replace the notion of Orlicz functions in
our discussions.

\Definition: A {\dt \phifunction\/} is a function $F:[0,\infty) \to
[0,\infty)$\ such that
\itemi $F(0) = 0$;
\itemii $\lim_{n\to\infty} F(t) = \infty$;
\itemiii $F$\ is strictly increasing;
\itemiv $F$\ is continuous;
\moreproclaim\noindent
We will say that a \phifunction\ $F$\ is {\dt dilatory\/} if for some
$c_1,c_2 > 1$\ we have $F(c_1 t)\ge c_2 F(t)$\ for all $0\le t<\infty$. We
will say that $F$\ satisfies the {\dt \Deltacond\/} if $F^{-1}$\ is dilatory.
\moreproclaim
If $F$\ is a \phifunction, we will define the function $\tilde F(t)$\ to
be $1/F(1/t)$\ if $t>0$, and $0$\ if $t=0$.

The definition of a \phifunction\ is slightly more restrictive than that of an
Orlicz function in that we insist that $F$\ be strictly increasing.
The notion of dilatory replaces
the notion of convexity. The Orlicz spaces generated by dilatory
functions are only quasi-Banach spaces, in contrast to those generated by Orlicz
functions, which are Banach spaces. 
The \Deltacond\ appears widely in literature
about Orlicz spaces.

\Definition: If $(\Omega,\F,\mu)$\ is a measure space, and $F$\ is a
\phifunction, then we define the {\dt Luxemburg functional\/} by
$$ \normo{f}_F = \inf\left\{\, c :
   \int_\Omega F\bigl(\modo{f(\omega)}/c\bigr) \,d\mu(\omega) \le 1 \,\right\} 
   ,$$
for every measurable function $f$. We define the {\dt Orlicz space},
$L_F(\Omega,\F,\mu)$\ (or $L_F(\mu)$, $L_F(\Omega)$ or $L_F$\ for short), to be
the vector space of measurable functions $f$\ for which $\normo f_F < \infty$,
modulo functions that are zero almost everywhere.

Now we define the Orlicz--Lorentz spaces.

\Definition: If $(\Omega,\F,\mu)$\ is a measure space, and
$F$\ and $G$\ are \phifunction s, then we define the {\dt
Orlicz--Lorentz functional\/} of a measurable function $f$\ by 
$$ \normo{f}_{F,G} = \normo{f^*\circ\tilde F\circ\tilde G^{-1}}_G .$$
We define the {\dt Orlicz--Lorentz space},
$L_{F,G}(\Omega,\F,\mu)$\ (or $L_{F,G}(\mu)$, $L_{F,G}(\Omega)$ or $L_{F,G}$\ for
short), to be the vector space of measurable functions $f$\ for which $\normo
f_{F,G} <\infty$, modulo functions that are zero almost everywhere.

\Definition: If $(\Omega,\F,\mu)$\ is a measure space, and
$F$\ is a \phifunction, then we define the ({\dt weak-}){\dt
Orlicz--Lorentz functional\/} by 
$$ \normo{f}_{F,\infty} = \sup_{x\ge0} \tilde F^{-1}(x) f^*(x) .$$
We define the {\dt Orlicz--Lorentz space},
$L_{F,\infty}(\Omega,\F,\mu)$\ (or $L_{F,\infty}(\mu)$, $L_{F,\infty}(\Omega)$ or
$L_{F,\infty}$\ for short), to be the vector space of measurable functions $f$\
for which $\normo f_{F,\infty} <\infty$, modulo functions that are zero almost
everywhere.

We see that $L_{F,F} = L_F$\ with equality of norms, and that if $F(t) =
t^p$\ and $G(t) = t^q$, then $L_{F,G} = L_{p,q}$, and $L_{F,\infty} =
L_{p,\infty}$, also with equality of norms. Thus, if $F(t) = t^p$, we shall
write $L_{p,G}$\ for $L_{F,G}$, and $L_{G,p}$\ for $L_{G,F}$.
Also, if $A$\ is any measurable set, then $\normo{\chi_A}_{F,G} =
\normo{\chi_A}_{F,\infty} = \normo{\chi_A}_F = \tilde F^{-1}\bigl(\mu(A)\bigr)$.

The Orlicz--Lorentz spaces defined here are equivalent to the
definition given in the introduction, as we now describe.

\Definition: A {\dt weight function\/} is a function $w:(0,\infty) \to
(0,\infty)$\ such that
$$ W(t) = \int_0^t w(s)\,ds  $$
is a \phifunction.

Then if $w$\ is a weight function, and $G$\ is a \phifunction, then 
$\Lambda_{w,G} = L_{\tilde W^{-1} \circ G,G}$, where
$$ W(t) = \int_0^t w(s) \, ds .$$

Now let us provide some examples. We define the {\dt modified
logarithm\/} and the {\dt modified exponential\/} functions by
$$ \eqalignno{
   \lm(t) 
   &=\cases{ 1+\log t                  & if $t\ge1$\cr
             1/\bigl(1+\log(1/t)\bigr) & if $0<t<1$\cr
             0                         & if $t=0$;\cr}\cr
   \em(t) = \lm^{-1}(t)
   &=\cases{ \exp(t-1)                 & if $t\ge1$\cr
             \exp\bigl(1-(1/t)\bigr)   & if $0<t<1$\cr
             0                         & if $t=0$.\cr}\cr}$$
These functions are designed so that for large $t$\ they behave like
the logarithm and the exponential functions, so that $\lm 1=1$\ and $\em 1=1$,
and so that $\widetilde{\lm}=\lm$\ and $\widetilde{\em} = \em$. Then the
functions $t^p(\lm t)^\alpha$\ and $\em(t^p)$\ are \phifunction s whenever
$0<p<\infty$\ and $-\infty<\alpha<\infty$. If the measure space is a
probability space, then the Orlicz spaces created using these functions are also
known as {\dt Zygmund spaces}, and the Orlicz--Lorentz spaces $L_{t^p(\lm
t)^\alpha,q}$\ and $L_{\em(t^p),q}$\ are known as {\dt Lorentz--Zygmund
spaces} (see, for example, \ref[B--S]).

Finally, we define the notions of equivalence.

\Definition: We say that two \phifunction s $F$\ and $G$\ are {\dt
equivalent\/} (in symbols $F\asymp G$) if for some number $c<\infty$\ we have
that $F(\invc t) \le G(t) \le F(ct)$\ for all $0\le t<\infty$.
\moreproclaim

\Definition: We say that two function spaces $X$\ and $Y$\ on the same
measure space are {\dt equivalent\/} if for some number $c<\infty$\ we have that
$f\in X \Leftrightarrow f\in Y$\ with $\invc \normo f_X \le \normo f_Y \le c\,
\normo f_X$\ for all measurable functions $f$.

\beginsection 3.\enspace Survey of Known Comparison Results

There are at least four obvious questions about Orlicz--Lorentz spaces.
\itemi For which \phifunction s $F$\ and $G$\ is $L_{F,G}$\ equivalent to a
normed space (or $p$-convex, or $q$-concave)?
\itemii What are the Boyd indices of the Orlicz--Lorentz spaces?
\itemiii What are necessary and sufficient conditions for $L_{F_1,G_1}$\ and
$L_{F_2,G_2}$\ to be equivalent?
\itemiv Is every rearrangement space equivalent to some Orlicz--Lorentz space?

\noindent
The first and second questions are intimately related, and will be dealt with
in another paper \ref[Mo2]. In general, they are very hard to answer. The
third question is the subject of this paper. As a corollary,
we will also be able to answer the fourth question.

There have already been many comparison results for Lorentz spaces. Indeed,
Lorentz himself provided one of the first in 1961 \ref[Lo3]. He found
necessary and sufficient conditions for $\Lambda_{w,1}$\ to be equivalent to an
Orlicz space. 

\Definition: A weight function $w$\ is said to be {\dt strictly monotone\/} if
either
\itemi $w$\ is strictly increasing, $w(t)\to0$\ as $t\to0$\ and
$w(t)\to\infty$\ as $t\to\infty$,
\moreproclaim\noindent
or
\item{i${}'$)} $w$\ is strictly decreasing, $w(t)\to\infty$\ as $t\to0$\ and
$w(t)\to0$\ as $t\to\infty$.
\moreproclaim\noindent

\Definition: A strictly monotone weight function is said to satisfy
{\dt\conditionL\/}\ if there is a number $c<\infty$\ such that
$$ \int_0^\infty {dt\over w^{-1}\bigl(c\,w(t)\bigr) } < \infty .$$

\proclaim Theorem 3.1. Let $w:(0,\infty)\to(0,\infty)$\ be a decreasing,
strictly monotone weight function. Then the following are equivalent. 
\itemi $\Lambda_{w,1}$\ is equivalent to an
Orlicz space. 
\itemii $w$\ satisfies \conditionL.

The kinds of weight functions that satisfy \conditionL\ are slowly increasing or
slowly decreasing functions. An example that Lorentz implicitly gave is
$$ w(t) = \cases{ t^{-1/\log(1+\log t)} &if $t\ge 1$\cr
                  t^{1/\log(1-\log t)} &if $0<t<1$.\cr}$$

Recently, Raynaud \ref[R] noticed that the above result is also true if $w$\
is strictly increasing. He then went on to show the following result.

\proclaim Theorem 3.2. Let $w$\ be a weight function, and $0<p<\infty$. If
there are strictly monotone weight functions $w_0$\ and $w_1$\ satisfying
\conditionL\ and a number $c<\infty$\ such that
$$ \invc {W(t)\over t} \le w_0(t) w_1(t) \le c\, {W(t)\over t} ,$$
then $\Lambda_{w,p}$\ is equivalent to an Orlicz space.

It may seem that the scope of these
results is limited, but this is not really the case. Using Lemmas~5.1.2
and~5.1.3 below, one can use these results to find sufficient conditions for
equivalence of two Orlicz--Lorentz spaces that are no stronger than the
conditions given in this paper.

There are also results due to
Bennett and Rudnick \ref[B--R] (see also \ref[B--S]). They proved the following
results for probability spaces, but using their methods, it is not too hard to
see that these results are true for all measure spaces.

\proclaim Theorem 3.3. For every $0<p<\infty$, and for every
$-\infty<\alpha<\infty$, we have that $L_{t^p(\lm t)^\alpha}$\ and $L_{t^p(\lm
t)^\alpha,p}$\ are equivalent.

\proclaim Theorem 3.4. For every $\beta>0$, we have that $L_{\em(t^\beta)}$\ and
$L_{\em(t^\beta),\infty}$\ are equivalent.

\beginsection 4.\enspace Comparison of Orlicz--Lorentz Spaces

In this section, we state the main results of this paper, and give necessary and
sufficient conditions for which, given certain restrictions upon $G_1$\ and
$G_2$, we have that $\normo f_{F,G_1} \le c\, \normo f_{F,G_2}$. Thus we 
find necessary and sufficient conditions for $L_{F_1,G_1}$\ and $L_{F_2,G_2}$\
to be equivalent. 

We first notice that $\normo f_{p,q_1} \le \normo f_{p,q_2}$\ whenever $q_1 \ge
q_2$\ (see \ref[H]). This suggests that we have a result something like:\enspace
if $G_1\circ G_2^{-1}$\ is a convex function, then $\normo f_{F,G_1} \le c\,
\normo f_{F,G_2}$. And this is indeed the case. However, more is true.
For example, if $G(t) = t\lm t$, then it follows from Theorem~3.3 that
$L_{G,1}$\ is equivalent to $L_{G,G}$. Thus, it would seem that we only need to
know that $G_1\circ G_2^{-1}$\ is `close,' in some sense, to a convex function.

In this paper, we establish precisely what this notion of closeness is.
But, before stating the conditions, we first give a little bit of motivation. We
note that a dilatory \phifunction\ $G$\ is determined completely, up to
equivalence, by its values $G(a^n)$, where $a>1$\ is any fixed number, and $n$\
ranges over all integers. Thus, we note that a \phifunction\ $G$\ is
equivalent to a convex function if and only if for some $a>1$\ and $N\in\N$, and
all $n\in\Z$\ and $m\in\N$, we have that $G(a^{n+m}) \ge a^{m-N} G(a^n)$\ (see
Lemma~5.4.2 below).

In all that follows, we take the natural numbers to be $\N = \{1$, $2$,
$3,\ldots\}$.

\Definition: Let $G$\ be a \phifunction. We say that $G$\ is
\itemi {\dt almost convex\/} if there are numbers $a>1$, $b>1$\ and
$N\in\N$\ such that for all $m\in\N$, the cardinality of the set of $n\in\Z$\
such that we do {\bf not} have
$$ G(a^{n+m}) \ge a^{m-N} G(a^n) $$
is less than $b^m$;
\itemii {\dt almost concave\/} if there are numbers $a>1$, $b>1$\ and
$N\in\N$\ such that for all $m\in\N$, the cardinality of the set of $n\in\Z$\
such that we do {\bf not} have
$$ G(a^{n+m}) \le a^{m+N} G(a^n) $$
is less than $b^m$;
\itemiii {\dt almost linear\/} if there are numbers $a>1$, $b>1$\ and
$N\in\N$\ such that for all $m\in\N$, the cardinality of the set of $n\in\Z$\
such that we do {\bf not} have
$$ a^{m-N} G(a^n) \le G(a^{n+m}) \le a^{m+N} G(a^n) $$
is less than $b^m$;
\itemiv {\dt almost constant\/} if there are numbers $a>1$, $b>1$\ and
$N\in\N$\ such that for all $m\in\N$, the cardinality of the set of $n\in\Z$\
such that we do {\bf not} have
$$ G(a^{n+m}) \le a^N G(a^n) $$
is less than $b^m$;
\itemv {\dt almost vertical\/} if $G^{-1}$\ is almost constant.

We will also express our results in terms of what we shall call \conditionJ.

\Definition: If $F$\ and $G$\ are \phifunction s, then say that $F$\ is {\dt
equivalently less convex than\/} $G$\ (in symbols $F \prec G$) if $G\circ
F^{-1}$\ is equivalent to a convex function. We say that $F$\ is {\dt
equivalently more convex than\/} $G$\ (in symbols $F \succ G$) if $G$\ is
equivalently less convex than $F$.

\Definition:
A \phifunction\ $F$\ is said to be an {\dt \Nfunction\/}
if it is equivalent to a \phifunction\ $F_0$\ such that $F_0(t)/t$\ is
strictly increasing, $F_0(t)/t \to \infty$\ as $t\to \infty$, and $F_0(t)/t \to
0$\ as $t\to0$.

\Definition:
A \phifunction\ $F$\ is said to be {\dt complementary\/}
to a \phifunction\ $G$\ if for some $c<\infty$\ we have 
$$ \invc t \le F^{-1}(t) \cdot G^{-1}(t) \le ct 
   \qquad (0\le t<\infty) .$$
If $F$\ is an \Nfunction, we will let $F^*$\ denote a function complementary
to $F$.

The notation $F^*$\ makes sense if $F$\ is an \Nfunction, because then there
is always a function $G$\ complementary to $F$, and further, if $G_1$\ and $G_2$\
are both complementary to $F$, then $G_1$\ and $G_2$\ are equivalent.

Our definition of a complementary function differs from the usual definition. If
$F$\ is an \Nfunction\ that is convex, then the complementary function is
usually defined by $F^*(t) = \sup_{s\ge0} \bigl(st-F(s)\bigr)$. However, it is
known that $t\le F^{-1}(t) \cdot {F^{*}}^{-1}(t) \le 2t$\ (see \ref[K--R]). Thus
our definition is equivalent.

\Definition: An \Nfunction\ $H$\ is said to satisfy {\dt \conditionJ\/}\ if
$$ \normo{1/\tilde H^*{}^{-1}}_{H^*} < \infty .$$

The kinds of \Nfunction s that satisfy \conditionJ\ are slowly rising functions.
These are essentially the kinds of Orlicz functions that Lorentz describes in
Theorem~1 of his paper \ref[Lo3].

We also describe our results in a third fashion. The following definitions are
motivated by the fact that $G_1\succ G_2$\ if
and only if for some $c<\infty$\ and all $s\ge 1$\ and $t>0$\ we have that
$G_1(st)/G_1(t) \ge \invc G_2(st)/G_2(t) $\ (see Lemma~5.4.2).

\Definition: Let $G_1$\ and $G_2$\ be \phifunction s. We say that
\itemi $G_1$ {\dt is almost less convex than\/} $G_2$\ if there are numbers
$a>1$, $b>1$\ and $N\in\N$\ such that for all $m\in\N$, the cardinality of the
set of $n\in\Z$\ such that we do {\bf not} have
$$ {G_1(a^{n+m})\over G_1(a^n)} \le a^{N} {G_2(a^{n+m})\over G_2(a^n)} $$
is less than $b^m$;
\itemii $G_1$ {\dt is almost more convex than\/} $G_2$\ if $G_2$ is almost
less convex than $G_1$;
\itemiii $G_1$ {\dt is almost equivalent to\/} $G_2$\ if there are numbers $a>1$,
$b>1$\ and $N\in\N$\ such that for all $m\in\N$, the cardinality of the set of
$n\in\Z$\ such that we do {\bf not} have
$$ a^{-N} {G_2(a^{n+m})\over G_2(a^n)} \le
   {G_1(a^{n+m})\over G_1(a^n)} \le a^{N} {G_2(a^{n+m})\over G_2(a^n)} $$
is less than $b^m$.

Now we collect together the comparison results. For all these results, we
will assume that the measure space is $[0,\infty)$\ with Lebesgue measure.
In fact, any non-atomic infinite measure space will do. There are also similar
results for non-atomic probability spaces, and $\N$\ with the counting measure
(i.e. sequence spaces). We do not give details for these cases. However the
idea is that for non-atomic probability spaces, we need only consider the
properties of the relevant \phifunction s $G(t)$\ for large $t$, and for
sequence spaces, their properties for small $t$. Obviously, if one is only
interested in sufficient conditions for Orlicz--Lorentz spaces to be equivalent,
one can use any measure space. (Recall that $\asymp$\ means {\it
equivalent to}, see Section~2.)

\proclaim Proposition 4.1. Let $F_1$, $F_2$, $G_1$\ and $G_2$\ be 
\phifunction s.
\itemi If $F_1\circ G_1^{-1} \asymp F_2\circ G_2^{-1}$\ and $G_1\asymp G_2$,
and if one of $G_1$\ or $G_2$\ is dilatory, then $L_{F_1,G_1}$\ and
$L_{F_2,G_2}$\ are equivalent.
\itemii If $F_1\asymp F_2$, then $L_{F_1,\infty}$\ and $L_{F_2,\infty}$\ are
equivalent.
\moreproclaim

\proclaim Theorem 4.2. Let $F$, $G_1$\ and $G_2$\ be \phifunction s. Consider
the following statements.
\itemi For some $c<\infty$, we have that 
$\normo f_{F,G_1} \le c\, \normo f_{F,G_2}$\ for all measurable $f$.
\itemii $G_1\circ G_2^{-1}$\ is almost convex.
\itemiii There is an \Nfunction\ $H$\ satisfying \conditionJ\ such that $G_1\circ
G_2^{-1} \succ H^{-1}$. 
\itemiv $G_1$\ is almost more convex than $G_2$.
\moreproclaim\noindent
Then, if one of
$G_1$\ or $G_2$\ is dilatory, then we have (ii)\implies(i). If one of 
$G_1$\ or $G_2$\ is dilatory and $G_1$\ satisfies the \Deltacond, or if
$G_2$\ satisfies the \Deltacond, then we have (i)\implies(ii). If $G_2$\ is
dilatory and satisfies the \Deltacond, then we have (ii)\iff(iv). 
We always have (ii)\iff(iii).

\proclaim Theorem 4.3. Let $F_1$, $F_2$, $G_1$\ and $G_2$\ be \phifunction s
such that one of $G_1$\ or $G_2$\ is dilatory, and that one of $G_1$\ or $G_2$\
satisfies the \Deltacond. Then the following are equivalent.
\itemi $L_{F_1,G_1}$ and $L_{F_2,G_2}$\ are equivalent.
\itemii $F_1 \asymp F_2$, and $G_1\circ G_2^{-1}$\ is almost linear.
\itemiii $F_1 \asymp F_2$, and there exist \Nfunction s $H$\ and $K$\
satisfying \conditionJ\ such that $G_1\circ G_2^{-1} = H\circ K^{-1}$.
\itemiv $F_1 \asymp F_2$, and there exist \Nfunction s $H$\ and $K$\
satisfying \conditionJ\ such that $G_1\circ G_2^{-1} = H^{-1}\circ K$. 
\itemv $F_1 \asymp F_2$, and $G_1$\ is almost equivalent to $G_2$. 
\itemvi $F_1\asymp F_2$, and there exist \Nfunction s $H$\ and $K$\
satisfying \conditionJ\ and a number $c<\infty$\ such that  $ \invc G_1/G_2 \le
H/K \le c\, G_1/G_2 $.
\itemvii $F_1 \asymp F_2$, and there exist strictly monotone weight functions
$w_0$\ and $w_1$\ satisfying \conditionL\ and a number $c<\infty$\ such that
$\invc G_1/G_2 \le w_0 w_1 \le c\, G_1/G_2 $. 
\itemviii $F_1 \asymp F_2$, and there exists an almost linear \phifunction\ $F$\
and a number $c<\infty$\ such that
$\invc F(t) \le t G_1(t) / G_2(t) \le c\, F(t)$\ for all $t>0$.
\moreproclaim

The condition that one of the \phifunction s $G_1$\ or $G_2$\ satisfy the
\Deltacond\ is required, as is shown by the following example. Let $G_1(t) =
\em t$\ and $G_2(t) = \em t^2$. By Theorem~4.6 below, $L_{1,G_1}$\ and
$L_{1,G_2}$\ are both equivalent to $L_{1,\infty}$. But, it is clear that
$G_1\circ G_2^{-1}$\ is far from being almost linear. The author
doesn't know whether the condition
that one of $G_1$\ or $G_2$\ be dilatory is needed.

We are also able to obtain certain results stating that
in order to compare $L_{F_1,G_1}$\ and $L_{F_2,G_2}$, we need only compare the
norms for a certain class of test functions.

\Definition:
Let $\T_1$\ be the set of functions $f:[0,\infty)\to[0,\infty)$\ such that for
some $0=a_0<a_1<a_2<\ldots<a_n$\ we have that
$$ f(x) = \cases{
   1/a_i & if $a_{i-1}\le x < a_i$\ and $1\le i\le n$ \cr
   0     & otherwise. \cr }$$
If $F$\ is a \phifunction, let $\T_F = \{\, F^{-1}\circ f:f\in\T_1\} = \{\,
f\circ\tilde F^{-1}:f\in\T_1\}$.

\proclaim Theorem 4.4. Let $F$, $G_1$\ and $G_2$\ be \phifunction s. Suppose
that $G_2$\ is dilatory, and that one of $G_1$\ or $G_2$\
satisfies the \Deltacond. Then the following are equivalent.
\itemi For some $c<\infty$ we have that
$\normo f_{F,G_1} \le c\, \normo f_{F,G_2}$\ whenever $f^*\in\T_F$.
\itemii For some $c<\infty$ we have that
$\normo f_{F,G_1} \le c\, \normo f_{F,G_2}$\ for all measurable $f$.

\proclaim Theorem 4.5. Let $F_1$, $F_2$, $G_1$\ and $G_2$\ be \phifunction s.
Suppose that one of $G_1$\ or $G_2$\ is dilatory, and that one of $G_1$\ or
$G_2$\ satisfies the \Deltacond. Then the following are equivalent.
\itemi For some $c<\infty$ we have that
$\invc \normo f_{F_1,G_1} \le \normo f_{F_2,G_2} \le c\, \normo f_{F_1,G_1}$\
whenever $f^*\in\T_{F_1}$. 
\itemii For some $c<\infty$ we have that
$\invc \normo f_{F_1,G_1} \le \normo f_{F_2,G_2} \le c\, \normo f_{F_1,G_1}$\
whenever $f^*\in\T_{F_2}$. 
\itemiii $L_{F_1,G_1}$\ and $L_{F_2,G_2}$\ are
equivalent.
\moreproclaim

Finally, we give a result for the weak-Orlicz--Lorentz spaces.

\proclaim Theorem 4.6. Let $F_1$, $F_2$\ and $G$\ be \phifunction s. Then the
following are equivalent.
\itemi $L_{F_1,G}$\ and $L_{F_2,\infty}$\ are equivalent.
\itemii $F_1 \asymp F_2$, and $G$\ is almost vertical.
\itemiii $F_1 \asymp F_2$, and $\normo{1/\tilde G^{-1}}_G < \infty$.
\itemiv $F_1 \asymp F_2$, and $\tilde G^{-1}$\ satisfies \conditionL.
\moreproclaim

It is clear that all the results given in Section~3 follow from these results.
We are also able to answer a question of Raynaud, and prove the converse to
Theorem~3.2.

\proclaim Theorem 4.7. Let $w$\ be a weight function, and $0<p<\infty$. If
$\Lambda_{w,p}$\ is equivalent to an Orlicz space, then there are strictly
monotone weight functions $w_0$\ and $w_1$\ satisfying \conditionL\ and a number
$c<\infty$\ such that 
$$ \invc {W(t)\over t} \le w_0(t) w_1(t) \le c\, {W(t)\over t} .$$

\Proof: This follows immediately from the implication (i)\implies(vii) in
Theorem~4.3, and from the observation that a strictly monotone weight
function $w$\ satisfies \conditionL\ if and only if $w^p$\ satisfies
\conditionL\ for any $0<p<\infty$. 
\endproof

\beginsection 5.\enspace The Proof of the Results of Section 4

The proofs of the results of Section~4 are rather long. We will split the proof
into many lemmas that are grouped into several subsections
according to their nature. Many of the lemmas, if not obvious, are at least
`believable without proof,' and the reader may pass over them quickly. The key
results are contained in Sections~5.3, 5.5 and~5.6.

These proofs could be shortened considerably if we assumed throughout
that all \phifunction s were dilatory and satisfied the \Deltacond, but
then our results would be correspondingly weaker.  In particular, Theorem~6.1
below would be much less general.

\beginsection 5.1.\enspace The Elementary Propositions

The first result is obvious, and requires no proof.

\proclaim Lemma 5.1.1. Let $G$\ be a \phifunction.
\itemi If $G$\ is dilatory, then for all $c_1<\infty $\ there is a number
$c_2<\infty$\ such that if
$$ \int_\Omega G\circ f(\omega) \,d\mu(\omega) \le c_1 ,$$
then $\normo f_G \le c_2$.
\itemii If $G$\ is dilatory, then for all $c_1>0 $\ there is a number
$c_2>0$\ such that if
$$ \int_\Omega G\circ f(\omega) \,d\mu(\omega) \ge c_1 ,$$
then $\normo f_G \ge c_2$.

Now we have the first result from Section~4.

\Proofof Proposition 4.1: This is a simple consequence of Lemma~5.1.1.
\endproof

The following results describe the basic `algebra' that the Orlicz--Lorentz
spaces satisfy. Essentially, they allow one to reduce comparison of
Orlicz--Lorentz spaces to the problem of comparing $L_{1,G}$\ to $L_1$.
The proofs are straightforward, so we omit them.

\proclaim Lemma 5.1.2. Suppose that $F$, $G_1$\ and $G_2$\ are \phifunction s.
Then for any number $c<\infty$\ we have that $\normo f_{F,G_1} \le c\, \normo
f_{F,G_2}$\ for all measurable $f$\ (respectively, $f\in\T_F$) if and only if
$\normo f_{1,G_1} \le c\, \normo f_{1,G_2}$\ for all measurable
$f$\ (respectively, $f\in\T_1$).

\proclaim Lemma 5.1.3. Suppose that $G_1$, $G_2$\ and $H$\ are \phifunction s. 
\itemi If $H$\ is dilatory, then if for some number $c_1<\infty$\ we have that
$\normo f_{1,G_1} \le c_1 \normo f_{1,G_2}$\ for all measurable $f$\
(respectively, $f\in\T_1$), then for some number $c_2<\infty$\ we have that
$\normo f_{1,G_1\circ H} \le c_2 \normo f_{1,G_2\circ H}$\ for all measurable
$f$\ (respectively, $f\in\T_1$).
\itemii If $H$\ satisfies the \Deltacond, then if for some number
$c_1<\infty$\ we have that $\normo f_{1,G_1\circ H} \le c_1 \normo f_{1,
G_2\circ H}$\ for all measurable $f$\ (respectively, $f\in\T_1$), then for some
number $c_2<\infty$\ we have that $\normo f_{1,G_1} \le c_2 \normo f_{1,G_2}$\
for all measurable $f$\ (respectively, $f\in\T_1$).
\moreproclaim

\beginsection 5.2.\enspace Conditions for Functions to be Dilatory, etc

Here we collect the results that pertain to when a \phifunction\ is dilatory or
satisfies the \Deltacond. The first result is obvious.

\proclaim Lemma 5.2.1. Let $G$\ be a \phifunction.
\itemi If there are numbers $a>1$, $c_1>1$\ and $c_2 > 1$\ such that 
$$ c_1 G(a^n) \le G(c_2 a^n) $$
except for finitely many $n$, then $G$\ is dilatory.
\itemii If there are numbers $a>1$, $c_1>1$\ and $c_2 > 1$\ such that 
$$ c_1 G(a^n) \ge G(c_2 a^n) $$
except for finitely many $n$, then $G$\ satisfies the \Deltacond.
\moreproclaim

Now we show how the property of $G$\ being dilatory or satisfying the \Deltacond\
may be captured by the properties of $L_{F,G}$.

\proclaim Lemma 5.2.2. Suppose that $G$\ is a \phifunction. Then the following
are equivalent.
\itemi $G$\ is dilatory.
\itemii There is a number $c<\infty$\ such that we have $\normo f_{1,G} \le
c$\ for all functions $f:\R\to\R$\ of the form
$$ f(x) = \cases{ a^{-1} & if $0\le x < a$ \cr
                  b^{-1} & if $a\le x < b$ \cr
                  0      & otherwise, \cr}$$
where $b>a>0$.
\moreproclaim

\Proof:
First we will show that (i)\implies(ii). Given a function of the above form, we
note that
$$ \int_0^\infty G\circ f^*\circ \tilde G^{-1}(x) \,dx
   \le \tilde G(a) G(a^{-1})
     + \tilde G(b) G(b^{-1})
   = 2 .$$
Then the result follows immediately from Lemma~5.1.1.

To show that (ii)\implies(i), we will consider functions of the form
$$ f(x) = \cases{ G^{-1}(3^m) & if $0\le x < 1/G^{-1}(3^m)$ \cr
                  G^{-1}(3^n) & if $1/G^{-1}(3^m) \le x < 1/G^{-1}(3^n)$ \cr
                  0      & otherwise ,\cr}$$
where $m>n$\ are integers. Then we know that $\normo f_{1,G} \le c$, and so
$$ 1 \ge 
   \int_0^\infty G\bigl(\invc f^*\circ\tilde G^{-1}(x)\bigr) \,dx
   \ge 3^{-m} G\bigl(\invc G^{-1}(3^m)\bigr)
   + \ts{2\over 3}\, 3^{-n} G\bigl(\invc G^{-1}(3^n)\bigr) .$$
Therefore, for all except one $n\in\N$\ we have that
$$ \ts{2\over3} 3^{-n} G\bigl(\invc G^{-1}(3^n)\bigr) \le \ts{1\over2} ,$$
that is,
$$ \invc G^{-1}(3^n) \le G^{-1}\bigl(\ts{3\over4} 3^n\bigr) .$$
By Lemma~5.2.1, $G^{-1}$\ satisfies the \Deltacond, and hence $G$\
is dilatory.
\endproof

\proclaim Lemma 5.2.3. If $F$, $G_1$\ and $G_2$\ are \phifunction s such that
$G_1$\ is dilatory, and such that for some $c<\infty$\ we have $\normo f_{F,G_2}
\le c\, \normo f_{F,G_1}$\ for all $f\in\T_F$, then $G_2$\ is dilatory.

\Proof: This follows immediately from Lemmas~5.1.2 and~5.2.2.
\endproof

\proclaim Lemma 5.2.4. Suppose that $G$\ is a \phifunction. Consider the
following statements.
\itemi $G$\ satisfies the \Deltacond.
\itemii Given $c>1$, there are numbers $d>1$\ and $N\in\N$ such that we have
$\normo f_{1,G} \ge c$\ for all functions $f:\R\to\R$\ of the form
$$ f(x) = \cases{ d^{-k_i} & if $d^{k_{i-1}}\le x < d^{k_i}$ and $1<i\le N$ \cr
                  0      & otherwise, \cr}$$
where $k_1<k_2<\ldots<k_N$\ are integers, and $k_0 = -\infty$.
\moreproclaim
Then (ii)\implies(i).  Furthermore, if $G$\ is dilatory, then (i)\implies
(ii).

\Proof:
First we will show that (i)\implies(ii) when $G$\ is dilatory.
Choose $d$\ so that $\tilde G(d^n) \ge
2 \tilde G(d^{n-1})$\ for all $n\in\Z$. Then if $f$\ is of the above form,
we have that
$$ \int_0^\infty G\circ f^*\circ \tilde G^{-1}(x) \,dx
   \ge \sum_{i=1}^N \bigl(\tilde G(d^{k_i})-\tilde G(d^{k_i-1})\bigr) 
       G(d^{-k_i})
   \ge N/2 .$$
Thus if $G$\ satisfies the \Deltacond, then by Lemma~5.1.1, there is some
$N\in\N$\ such that for all $f$\ of the above form, we have that
$\normo f_{1,G} \ge c$.

To show that (ii)\implies(i), let us pick $c>2$. Then for any function of the
above form, we have that
$$ 1 < 
   \int_0^\infty G\bigl(2^{-1} f^*\circ\tilde G^{-1}(x)\bigr) \,dx
   \le \sum_{i=1}^N \tilde G(d^{k_i}) G(2^{-1} d^{-k_i}). $$
Therefore, the cardinality of the set of $n\in\N$\ such that
$$ \tilde G(d^n) G(2^{-1} d^{-n}) \le N^{-1} $$
is less than $N$. By Lemma~5.2.1, this shows that $G$\ satisfies the
\Deltacond.
\endproof

\proclaim Lemma 5.2.5. Suppose that $F$, $G_1$\ and $G_2$\ are \phifunction s 
such that one of $G_1$\ and $G_2$\ is dilatory.  If
$G_1$\ satisfies the \Deltacond, and for some $c<\infty$\ we have that
$\normo f_{F,G_2} \ge \invc \normo f_{F,G_1}$\ for all $f\in\T_F$, then
$G_2$\ satisfies the \Deltacond.

\Proof: By Lemma~5.2.3, we have that $G_1$\ is dilatory.  Now the 
result follows immediately from Lemmas~5.1.2 and~5.2.4.
\endproof

\beginsection 5.3.\enspace Comparison Conditions for $L_{1,G}$

In this subsection, we give the key lemma that demonstrates the relationship
between the almost convexity of $G$, and the comparison between $L_{1,G}$\ and
$L_1$. As a corollary, we will also obtain results that show that in the
definition of the `almost properties' that we can take the value of $a$\ to be
arbitrarily large.

\proclaim Lemma 5.3.1. Suppose that $G$\ is a \phifunction. Then
the following are equivalent.
\itemi For some $c<\infty$, we have that $\normo f_{1,G} \le c\, \normo f_1 $\
for all measurable $f$.
\itemii For some $c<\infty$, we have that $\normo f_{1,G} \le c\, \normo f_1 $\
for all $f\in\T_1$.
\itemiii For all sufficiently large $a$, there are numbers $b>1$\ and
$N\in\N$\ such that for all $m\in\N$, the cardinality of the set of $n\in\Z$\
such that we do {\bf not} have
$$ G(a^{n+m}) \ge a^{m-N} G(a^n) $$
is less than $b^m$.
\itemiv $G$\ is almost convex.

The proof will require the next lemma.

\proclaim Lemma 5.3.2. Let $G$\ be a \phifunction. If $G$\ is almost convex, then
given $a'>1$,
there are numbers $a>1$, $b>1$, $c<\infty$\ and $N\in\N$\ such that
$a>\max\{a',b\}$\ and such that for all $m\in\N$, the cardinality of the set
of $n\in\Z$\ such that we do {\bf not} have 
$$ G(a^{n+m}) \ge a^{m-N} G(a^n) $$
is less than $cb^m$.
\moreproclaim
There are similar results if $G$\ is almost concave or almost constant.

\Proof: There are numbers $a>1$, $b>1$\ and $N\in\N$\ such that for all $m\in\N$,
the cardinality of the set
$$ A_m = \{\, n\in\Z : G(a^{n+m}) < a^{m-N} G(a^n) \} $$
is less than $b^m$. Pick $c\in\N$\ such that $a^c > b$\ and $a^c>a'$, and let 
$$ A'_m = \{\, n\in\Z : G(a^{c(n+m)}) < a^{c(m-N)} G(a^{cn}) \} .$$
Then, if $n\in A'_m$, then at least one of $cn$, $cn+m,\ldots,$\ or $cn+(c-1)m
\in A_m$, and hence $\modo{A'_m} < cb^m$.
\endproof

\Proofof Lemma 5.3.1: Clearly, (i)\implies(ii) and (iii)\implies(iv). We will
show (ii)\implies(iii). By Lemma~5.2.3, we know that $G$\ is dilatory. Thus we
suppose that $a>2$\ and $G(at) \ge 2G(t)$\ for all $t\ge0$. Choose $N$\
so that $a^{N-1}>c$. We will prove the result by showing that there cannot be
numbers $m\in\N$\ and $n_1<n_2<\ldots<n_{a^{m-N+1}}$\ such that 
$$ G(a^{n_i}) < a^{m-N} G(a^{n_i-m}) .$$
For otherwise, consider the
function
$$ f(x) = \cases{ a^{n_i} & if $a^{-n_{i+1}} \le x < a^{-n_i}$\ 
                             and $1\le i \le a^{m-N+1}$ \cr
                  0        & otherwise, \cr}$$
where we take $n_{a^{m-N+2}} = \infty$.
Clearly, $\normo f_1 \le a^{m-N+1}$. But also, we have the following
inequalities. 
$$ \eqalignno{
   & \int_0^\infty G\bigl(a^{-m} f^*\circ \tilde G^{-1}(x)\bigr) \,dx \cr
   &\ge \sum_{i=1}^{a^{m-N+1}}
        \int_{\tilde G(a^{-n_{i+1}})}^{\tilde G(a^{-n_i})}
        G\bigl(a^{-m} f^*\circ \tilde G^{-1}(x)\bigr) \,dx \cr
   &\ge \half \sum_{i=1}^{a^{m-N+1}}
        {G(a^{-m+n_i})\over G(a^{n_i})} \cr
   &\ge a/2 > 1 ,\cr}$$
where the penultimate inequality follows because $\tilde G(a^{-n_i}) > 2
\tilde G(a^{-n_{i+1}})$.
Thus $\normo f_{1,G} \ge a^m \ge a^{N-1} \normo f_1$,
which is a contradiction.

Now we show that (iv)\implies(i).
By Lemma~5.3.2, there are numbers $a>b>1$, $c<\infty$\ and $N\in\N$\ such that
for all $m\in\N$, the cardinality of the set
$$ A_m = \{\, n\in\Z : G(a^n) < a^{m-N} G(a^{n-m}) \} $$
is less than $cb^m$.
Let $\{\lists k\}$\ be the (possible finite or
empty) set of integers not in $\bigcup_{m=1}^\infty A_m$. Define
the sequence of sets $B_m$\ by setting $B_1 = A_1$, and 
$$ B_m = \{k_{m-1}\} \cup A_m \setminus \bigcup_{m'<m} A_{m'} $$
for $m>1$. Then $\modo{B_m} \le cb^m$.

Choose $c_1 = (a-1)/a^{4}$. Suppose $f$\ is a measurable function such that 
$$ \normo f_1 = \int_0^\infty f^*(x) \,dx \le c_1 .$$
For each $n\in\Z$, let $m_n \in \Z\cup\{\infty\}$\ be such that 
$$ a^{-m_n} \le a^{-n} f^*(a^{-n}) \le a^{1-m_n} .$$
Then
$$ \eqalignno{
   c_1 
   &\ge \sum_{n=-\infty}^\infty \int_{a^{-n-1}}^{a^{-n}} f^*(x) \,dx \cr
   &\ge \sum_{n=-\infty}^\infty (a^{-n}-a^{-n-1}) f^*(a^{n}) \cr
   &\ge {a-1\over a} \sum_{n=-\infty}^\infty a^{-m_n} .\cr}$$
Therefore, 
$$ \sum_{n=-\infty}^\infty a^{-m_n} \le {a\over a-1} \, c_1 = a^{-3} .$$
In particular, we note that $m_n\ge 3$\ for all $n\in\Z$.

Then
$$ \eqalignno{
   &\int_0^\infty G\circ f^*\circ\tilde G^{-1}(x) \,dx \cr
   &= \sum_{n=-\infty}^\infty
      \int_{\tilde G(a^{-n-1})}^{\tilde G(a^{-n})}
      G\circ f^*\circ\tilde G^{-1}(x) \,dx \cr
   &\le \sum_{n=-\infty}^\infty 
        \tilde G(a^{-n}) G\bigl(f^*(a^{-n-1})\bigr) \cr
   &\le \sum_{n=-\infty}^\infty
        G(a^{2+n-m_{n+1}}) / G(a^{n}) .\cr}$$
Now let 
$$ V = \{\, n\in\Z : n\not\in B_m \hbox{ for all } m\le m_{n+1}-2 \} .$$
Then
$$ \eqalignno{
   &\int_0^\infty G\circ f^*\circ\tilde G^{-1}(x) \,dx \cr
   &\le \sum_{n\in V}
      G(a^{2+n-m_{n+1}}) / G(a^{n}) 
      + \sum_{m=1}^\infty \sum_{n\in B_m\setminus V}
             G(a^{2+n-m_{n+1}}) / G(a^{n}) .\cr} $$
If $n\in V$, then either $n\not\in A_{m_{n+1}-2}$, or $m_{n+1}=\infty$, and
so
$$ G(a^{2+n-m_{n+1}}) / G(a^{n}) \le a^{N+2-m_{n+1}} .$$
If $n\in B_m \setminus V$, then $m\le m_{n+1}-2$, and so
$$ G(a^{2+n-m_{n+1}}) / G(a^{n}) \le G(a^{n-m}) / G(a^{n}) .$$
If we also know that $m>1$, then $n\not\in A_{m-1}$, and so
$$ G(a^{2+n-m_{n+1}}) / G(a^{n}) \le G(a^{1+n-m}) / G(a^{n}) 
   \le a^{N+1-m} .$$
Therefore,
$$ \eqalignno{
   &\int_0^\infty G\circ f^*\circ\tilde G^{-1}(x) \,dx \cr
   &\le \sum_{n\in V}
      a^{N+2-m_{n+1}} 
      + \sum_{n\in B_1\setminus V}
             G(a^{n-1}) / G(a^{n}) 
      + \sum_{m=2}^\infty \sum_{n\in B_m\setminus V}
             a^{N+1-m} \cr
   &\le a^{N-1}
      + \sum_{n\in B_1}
             G(a^{n-1}) / G(a^{n}) 
      + \sum_{m=2}^\infty c b^m a^{N+1-m} ,\cr }$$
which is a finite number whose value does not depend on $f$. However, by
Lemma~5.2.1, $G$\ is dilatory, and hence by Lemma~5.1.1, $\normo f_{1,G}$\ is
bounded by some number that does not depend on $f$. 
\endproof

\beginsection 5.4.\enspace Convexity and Concavity Conditions

In this subsection, we give basic results about convexity and concavity, and
their `almost' equivalents. First, we give a technical lemma whose proof is
obvious.

\proclaim Lemma 5.4.1. Let $G$\ be a \phifunction. Define a map $f:\Z\to\Z$\
so that for all $n\in\Z$\ we have
$$ a^{f(n)-1} \le G(a^n) < a^{f(n)} .$$
\itemi If $G$\ is dilatory, then there is a number $L\in\N$\ such that
$f(n_1) \ne f(n_2)$\ if $\modo{n_1-n_2} \ge L$.
\itemii If $G$\ satisfies the \Deltacond, then there is a number $M\in\N$\ such
that 
$$ \modo{f(n_1)-f(n_2)} \le M\modo{n_1-n_2} $$
for all $n_1$, $n_2\in\Z$.
\moreproclaim

Next we give some results about convexity.

\proclaim Lemma 5.4.2. Let $G_1$\ and $G_2$\ be \phifunction s. Consider the
following statements.
\itemi $G_1$\ is equivalently more convex than $G_2$.
\itemii There is a number $c<\infty$\ such that $G_1\circ G_2^{-1}(st) \ge
\invc s G_1\circ G_2^{-1}(t) $\ for all $s\ge1$\ and $t\ge0$.
\itemiii There is a number $c<\infty$\ such that $G_1(uv)/G_1(v) \ge \invc
G_2(uv)/G_2(v)$\ for all $u\ge1$\ and $v>0$.
\itemiv There are numbers $a>1$\ and $N\in\N$\ such that for all $m\in\N$\ and
$n\in\Z$, we have that
$$ {G_1(a^{n+m}) \over G_1(a^n)} \ge a^{-N} {G_2(a^{n+m}) \over G_2(a^n)} .$$
\moreproclaim\noindent
Then we have (i)\iff(ii)\iff(iii)\implies(iv). If one of $G_1$\ or $G_2$\
satisfies the \Deltacond, then (iv)\implies(iii).

\Proof: The implications (i)\implies(ii) and (iii)\implies(iv) are obvious. The
implications (ii)\iff(iii) follow by setting $t=G_2(v)$\ and $st = G_2(uv)$. 

To
show (ii)\implies(i), we let 
$$ H_0(t) = \inf_{s>1} {G_1\circ G_2^{-1}(st) \over s} 
   \quad\hbox{and}\quad
   H(t) = \int_0^t {H_0(s)\over s} \,ds .$$
Then it is easy to see that $H$\ is convex, and that
$H$\ is equivalent to $G_1\circ G_2^{-1}$.

Now suppose that $G_1$\ 
satisfies the \Deltacond. We show (iv)\implies(iii).
Let $L\in\N$\ be such that $G_1(at) \le a^L G_1(t)$\ for all
$t\ge0$. Suppose that for some $u>1$\ and $v>0$\ we have that
$$ {G_1(uv) \over G_1(v)} < a^{-N-3L} {G_2(uv) \over G_2(v)} .$$
Let $m$\ and $n$\ be such that
$ a^m \le u < a^{m+1} $\ and $ a^n \le v < a^{n+1} $. Then
$$ a^{-3L} {G_1(a^{m+n+2}) \over G_1(a^n)}
   \le
   {G_1(a^{m+n}) \over G_1(a^{n+1})}
   <
   a^{-N-3L} {G_2(a^{m+n+2}) \over G_2(a^n)} ,$$
which is a contradiction.

The argument for when $G_2$\ satisfies the \Deltacond\ is similar.
\endproof

Now we start looking at the `almost properties.' First we relate almost
convexity to almost concavity.

\proclaim Lemma 5.4.3. Suppose that $G$\ is a \phifunction.
\itemi If $G$\ is almost convex and satisfies the \Deltacond, 
then $G^{-1}$\ is almost concave.
\itemii If $G$\ is almost concave, then $G^{-1}$\ is almost convex.
\moreproclaim

This will follow from the next lemma.

\proclaim Lemma 5.4.4. Suppose that $G$\ is a \phifunction\ and that $a>1$.
Consider the following statements.
\itemi There are numbers $b>1$\ and
$N\in\N$\ such that for all $m\in\N$, the cardinality of the set of $n\in\Z$\
such that we do {\bf not} have
$$ G(a^{n+m}) \ge a^{m-N} G(a^n) $$
is less than $b^m$.
\itemii There are numbers $b>1$\ and
$N\in\N$\ such that for all $m\in\N$, the cardinality of the set of $n\in\Z$\
such that we do {\bf not} have
$$ G^{-1}(a^{n+m}) \le a^{m+N} G^{-1}(a^n) $$
is less than $b^m$.
\moreproclaim\noindent
Then (ii)\implies(i). If, in addition, $G$\ satisfies the \Deltacond, then
(i)\implies(ii).

\Proof: We will show that (i)\implies(ii) when $G$\ satisfies the \Deltacond.
Let $f:\Z\to\Z$\ be defined so that 
$$ a^{f(n)-1} \le G^{-1}(a^n) < a^{f(n)} .$$
Since $G^{-1}$\ is dilatory, by Lemma~5.4.1, we know that there is a number $L$\
such that for every $n\in\Z$\ we have that $\modo{f^{-1}(\{n\})} \le L$.

Let 
$$ A_m = \{\, n: G^{-1}(a^{n+m}) > a^{m+1+N} G^{-1}(a^n) \}.$$
Then it can easily be shown that 
$$ f(A_m) \subseteq \{\, n : G(a^{n+m+N}) < a^m G(a^{n})\} ,$$
and hence $\modo{A_m} \le L b^{m+N} \le b_0^m$, where $b_0 = L b^{N+1}$.

To show (ii)\implies(i) is similar. Let $g:\Z\to\Z$\ be defined so that
$$ a^{g(n)} \le G(a^n) < a^{g(n)+1} .$$
Since $G$\ is almost convex, it follows that $G$\ is dilatory. 
Now the proof proceeds as in (i)\implies(ii).
\endproof

Next, we deal with the composition of the `almost' properties.  One of
the main problems here is that given two \phifunction s, each with an
`almost' property, is that the $a$\ from the definition of the `almost'
property for each \phifunction\ could be different.  Fortunately, 
we have already developed the tools to deal
with this.  First, for `almost convexity,' the implication (iv)\implies(iii)
in Lemma~5.3.1 tells us that the $a$\ may be any arbitrarily large
number.  If the \phifunction\ is dilatory, then Lemma~5.4.4 also
allows the $a$\ to be any arbitrarily large number for the `almost concavity'
property.  Finally, for other `almost' properties. Lemma~5.3.2 allows us to
choose the $a$\ to be larger than any given number.

Thus we have the following result.

\proclaim Lemma 5.4.5. Let $G$\ be a \phifunction. If $G$\ is almost convex and
almost concave, then $G$\ is almost linear.

\proclaim Lemma 5.4.6. Let $G_1$\ and $G_2$\ be \phifunction s.
\itemi If $G_1$\ and $G_2$\ are almost convex, then $G_1\circ G_2$\ is almost
convex.
\itemii If $G_1$\ and $G_2$\ are almost concave, and if $G_2$\ is dilatory,
then $G_1\circ G_2$\ is almost concave.

\Proof: First we will prove part~(ii). By the explanation given above, we
may suppose that for one $a>1$, there are numbers $b_1>1$, $N_1\in\N$,
$b_2>1$\ and $N_2\in\N$\ such that for all $m\in\N$,
the cardinality of the set of $n\in\Z$\ such that
$$ G_1(a^{n+m}) > a^{m+N_1} G_1(a^n) $$
is less than $b_1^m$, and the cardinality of the set of $n\in\Z$\ such that
$$ G_2(a^{n+m}) > a^{m+N_2} G_2(a^n) $$
is less than $b_2^m$.

Define a function $f:\Z\to\Z$\ so that for all $n\in\Z$\ we have
$$ a^{f(n)-1} \le G_2(a^n) < a^{f(n)} .$$
Then by Lemma~5.4.1, there is a number $L\in\N$\ such that $\modo{f^{-1}(\{n\})}
\le L$\ for all $n\in\Z$. 

Then we see that for all
$m\in\N$, the cardinality of the set of $n\in\Z$\ such that
$$ G_1\circ G_2(a^{n+m}) > a^{m+N_1+N_2+1} G_1\circ G_2(a^n) $$
is less than $b_1^{m+N_1+N_2+1} + L b_2^{m+N_2+1}$.
For, if the above holds, and
$$ G_2(a^{n+m}) \le a^{m+N_2} G_2(a^n) ,$$
then
$$ G_1(a^{f(n)-1+m+N_2+1}) > a^{m+N_2+1+N_1} G_1(a^{f(n)-1}) .$$
The result follows.

To show part~(i), we note that as $G_2$\ is almost convex, that we already know
that $G_2$\ is dilatory. Now the argument follows as in part~(i).
\endproof

Now, we prove two lemmas that are `almost' analogues of Lemma~5.4.2.

\proclaim Lemma 5.4.7. Let $G_1$\ and $G_2$\ be \phifunction s such that $G_2$\
is dilatory and satisfies the \Deltacond. Then we have the following.
\itemi $G_1\circ G_2^{-1}$\ is almost convex if and only if $G_1$\ is almost
more convex than $G_2$.
\itemii $G_1\circ G_2^{-1}$\ is almost concave if and only if $G_1$\ is almost
less concave than $G_2$.

\Proof: We will show that if $G_1\circ G_2^{-1}$\ is almost convex, then $G_1$\
is almost more convex than $G_2$. All the other assertions follow similarly.

So, there are numbers $a>1$, $b>2$\ and $N\in\N$\ such that for
all $m\in\N$, the cardinality of
$$ A_m = \{\, n :
   G_1\circ G_2^{-1}(a^{n+m}) < a^{m-N} G_1\circ G_2^{-1}(a^n) \} $$
is less than $b^m$. Let us define a map $f:\Z\to\Z$\ so that for all $n\in\Z$\
we have that
$$ a^{f(n)-1} \le G_2(a^n) < a^{f(n)} .$$
Then, by Lemma~5.4.1, there are numbers $L$, $M\in\N$\ such that
$f(n+L) > f(n)$\ for all $n\in\Z$, and such that $f(m+n) - f(n) \le
Mm$\ for all $m\in\N$\ and $n\in\Z$.

Now, for each $m\in\N$, let us consider the cardinality of the set 
$$ B_m = \left\{\, n:
   {G_1(a^{L(n+m)}) \over G_1(a^{Ln})} < a^{-N-2}
   {G_2(a^{L(n+m)}) \over G_2(a^{Ln})} \right\} .$$
If $n\in B_m$, let $n'=f(Ln)$\ and $m'=f\bigl(L(m+n)\bigr)-n'$. Then
$$ G_1\circ G_2^{-1}(a^{n'+m'-1}) < a^{m'-1-N} G_1\circ G_2^{-1}(a^{n'}) .$$
Clearly, this is impossible if $m'\le 1$, and otherwise, this implies that
$n'\in A_{m'-1}$. Since $m' \le Mm$, we see that
$$ \modo{B_m} < \sum_{m=1}^{Mm} b^{m'} \le b_0^m ,$$
where $b_0 = b^{M+1}/(b-1) $.
\endproof

\proclaim Lemma 5.4.8. Let $G_1$\ and $G_2$\ be \phifunction s such that
one of $G_1$\ or $G_2$\ is dilatory and one of $G_1$\ or $G_2$\ satisfies the
\Deltacond. Then $G_1\circ G_2^{-1}$\ is almost linear if and only if $G_1$\ is
almost equivalent to $G_2$.

\Proof: We note that in either case that if one of $G_1$\ or $G_2$\ is
dilatory, then both are, and if one of $G_1$\ or $G_2$\ satisfies the
\Deltacond, then both do. Now the proof proceeds as in Lemma~5.4.7.
\endproof

\beginsection 5.5.\enspace Condition ({\it L}) and Condition ({\it J})

In this subsection, we describe how the notions of satisfying \conditionL\
or \conditionJ\ relate to the `almost' properties.

\proclaim Lemma 5.5.1. Let $G$\ be a \phifunction. Then the following are
equivalent.
\itemi $G$\ is almost constant.
\itemii $\normo{1/\tilde G}_{G^{-1}} < \infty $.
\itemiii $\tilde G$\ satisfies \conditionL.
\moreproclaim

\Proof: This proof is very similar to the proof of Lemma~5.3.1, and so
we will omit many details.
First we show that (i)\implies(ii).
Following the same argument as the proof of (iv)\implies(i) in Lemma~5.3.1,
we construct numbers $a>b>1$, $c<\infty$\ and $N\in\N$\ and
a sequence of sets $B_m$\ such that $\modo{B_m} \le cb^m$, and such that
if $n\in B_m$\ for $m>1$ then
$$ G^{-1}\bigl(a^{-N}G(a^n)\bigr) \le a^{n-m+1} .$$
Hence,
$$ \eqalignno{
   & \int_0^\infty G^{-1}\bigl(a^{-N}/\tilde G(x)\bigr) \,dx \cr
   &= \sum_{n=-\infty}^\infty
      \int_{a^{-n}}^{a^{1-n}} G^{-1}\bigl(a^{-N}/\tilde G(x)\bigr) \,dx \cr
   &\le \sum_{m=1}^\infty \sum_{n\in B_m} a^{1-n} 
        G^{-1}\bigl(a^{-N} G(a^{n})\bigr) \cr
   &\le \sum_{n\in B_1} a^{1-n} 
        G^{-1}\bigl(a^{-N} G(a^{n})\bigr) 
      + \sum_{m=2}^\infty c b^m a^{1-n} a^{n-m+1} ,\cr }$$
which is a finite number. By Lemma~5.2.1, $G^{-1}$\ is dilatory, and so 
the result follows by Lemma~5.1.1.

That (ii)\implies(iii) is straightforward. To show that (iii)\implies(i), choose
$a>2$, and note that for some $N$, $M\in\N$\ we have that
$$ \int_0^\infty G^{-1}\bigl(a^{-N}/\tilde G(x)\bigr) \,dx \le a^M .$$
Then following a similar line of reasoning to that of the proof
of (ii)\implies(iii) in Lemma~5.3.1, it is possible to show that
there cannot be numbers $m\in\N$\ and $n_1<n_2<\ldots<n_{a^{m+2+M}}$\ such
that  
$$ G(a^{n_i}) > a^{N} G(a^{n_i-m}) .$$
\endproof

\proclaim Lemma 5.5.2. Let $H$\ is an \Nfunction. Then the following are
equivalent.
\itemi $H$\ is almost concave.
\itemii $H$\ satisfies \conditionJ.
\itemiii $\tilde H^*{}^{-1}$\ satisfies \conditionL.
\moreproclaim

\Proof: By Lemma~5.4.3, $H$\ is almost concave if and only if
$H^{-1}$\ is almost convex. Clearly, $H^{-1}$\ is almost convex if and only if
$H^*{}^{-1}$\ is almost constant. Now the result follows by Lemma~5.5.1.
\endproof

\beginsection 5.6.\enspace Condition ({\it J}) and the
`Almost' Properties

Now, we are ready to establish the relationship between being almost convex or
almost concave, and being more or less convex than some \Nfunction\ satisfying
\conditionJ.

\proclaim Lemma 5.6.1. Let $G$\ be a \phifunction. Then we have the following.
\itemi $G$\ is almost convex if and only if there is an \Nfunction\ $H$\
satisfying \conditionJ\ such that $G\succ H^{-1}$.
\itemii If $G$\ is almost concave, then there is an \Nfunction\ $H$\ satisfying
\conditionJ\ such that $G\prec H$.
\itemiii If there is an \Nfunction\ $H$\ satisfying \conditionJ\ such that
$G\prec H$, then $G^{-1}$\ is almost convex.
\moreproclaim

\Proof: We first note that if there is an \Nfunction\ $H$\ satisfying
\conditionJ\ such that either $G\succ H^{-1}$\ or $G^{-1}\prec H$, then by
Lemma~5.4.2, we have that $G$\ is almost convex. We will prove the other
implication of part~(i).

If $G$\ is almost convex, then we know that there are numbers $a>1$, $b>2$\ and
$N\in\N$\ such that for all $m\in\N$, the cardinality of the set
of $n\in\Z$\ such that
$$ G(a^{n+m}) < a^{m-N} G(a^n) $$
is less than $b^m$. Now we define a function $L:\{\,a^n:n\in\Z\} \to
(0,\infty)$\ by
$$ L(a^n) =
   \cases{\ds{ 
   \inf_{0=n_0<n_1<\ldots<n_K=n}\,
   \prod_{k=1}^K
   a^N \min\left\{ a^{n_k-n_{k-1}} ,\,
   {G(a^{n_k}) \over G(a^{n_{k-1}})} \right\}}
   &if $n\ge 0$ \cr
   \ds{
   \sup_{0=n_0>n_1>\ldots>n_K=n}\,
   \prod_{k=1}^K
   a^{-N} \max\left\{ a^{n_k-n_{k-1}} ,\,
   {G(a^{n_k}) \over G(a^{n_{k-1}})} \right\}}
   &if $n < 0$. \cr}$$
We may extend
the domain of $L$\ to $[0,\infty)$\ `log-linearly,' that is, by setting
$L(0)=0$, and
$$ L(a^nt) = L(a^n) \exp\left({\log t\over \log a}
   \log\left({L(a^{n+1}) \over L(a^n)}\right) \right) ,$$
for $n\in\Z$\ and $1\le t < a$.
We notice that $L(a^{n+1}) > L(a^n)$\ for all $n\in\Z$, and hence $L$\ is a
\phifunction.

Now, we note that if $m\in\N$\ and $n\in\Z$, then
$$ {L(a^{n+m})\over L(a^n)} \le
   \inf_{n=n_0<n_1<\ldots<n_K=n+m}\,
   \prod_{k=1}^K
   a^N \min\left\{ a^{n_k-n_{k-1}} ,\,
   {G(a^{n_k}) \over G(a^{n_{k-1}})} \right\} .$$
Thus, we have that $L(a^{n+m}) \le a^{m+N} L(a^n)$, and so, by Lemma~5.4.2,
$L^{-1}$\ is equivalent to a convex function. We also have that 
$$ {L(a^{n+m})\over L(a^n)} \le a^N {G(a^{n+m})\over G(a^n)} ,$$ 
and therefore, by Lemma~5.4.2 and since $L$\ satisfies the \Deltacond, $G \succ
L$.

We also notice that, since 
$$ \min\left\{ a^{n_2-n_0} ,\,
   {G(a^{n_2}) \over G(a^{n_0})} \right\}
   \ge
   \min\left\{ a^{n_2-n_1} ,\,
   {G(a^{n_2}) \over G(a^{n_1})} \right\}
   \min\left\{ a^{n_1-n_0} ,\,
   {G(a^{n_1}) \over G(a^{n_0})} \right\} $$
for $n_0 \le n_1 \le n_2$, we have that
$$ {L(a^{n+m})\over L(a^n)} \ge a^{-N}
   \inf_{n=n_0<n_1<\ldots<n_K=n+m}\,
   \prod_{k=1}^K
   a^N \min\left\{ a^{n_k-n_{k-1}} ,\,
   {G(a^{n_k}) \over G(a^{n_{k-1}})} \right\} .$$
Therefore, if $L(a^{n+m}) < a^{m-N} L(a^n)$, then for some $n\le n' < n'+m'
\le n+m$, we have that $G(a^{n'+m'}) < a^{m'-N} G(a^{n'})$. Therefore, the
cardinality of the set of $n\in\Z$\ satisfying $L(a^{n+m}) < a^{m-N} L(a^n)$\
is less than $m(b+b^2+\ldots+b^m)$, which is less than $b_0^m$\ for $b_0 =
2b^2/(b-1)$.

Therefore, $L$\ is almost convex. Now, we define the
\phifunction\ $H(t) = L^{-1}(t\lm t)$. It is clear that $L \succ H^{-1}$, and
hence $G\succ H^{-1}$. Since $t\lm t$\ is easily seen to be almost concave, it
follows by Lemmas~5.4.3 and~5.4.6 that $H$\ is almost convex. Clearly $H$\ is
an \Nfunction, and so by Lemma~5.5.2, we have that $H$\ satisfies
\conditionJ.

The proof of part~(ii) is similar. We know that there are numbers $a>1$, $b>2$\
and $N\in\N$\ such that for all $m\in\N$, the cardinality of the set
of $n\in\Z$\ such that
$$ G(a^{n+m}) > a^{m+N} G(a^n) $$
is less than $b^m$. Now we define a function $L:\{\,a^n:n\in\Z\} \to
(0,\infty)$\ by
$$ L(a^n) = 
   \cases{\ds{
   \sup_{0=n_0<n_1<\ldots<n_K=n}\,
   \prod_{k=1}^K
   a^{-N} \max\left\{ a^{n_k-n_{k-1}} ,\,
   {G(a^{n_k}) \over G(a^{n_{k-1}})} \right\}}
   &if $n\ge 0$ \cr
   \ds{
   \inf_{0=n_0>n_1>\ldots>n_K=n}\,
   \prod_{k=1}^K
   a^N \min\left\{ a^{n_k-n_{k-1}} ,\,
   {G(a^{n_k}) \over G(a^{n_{k-1}})} \right\}}
   &if $n < 0$, \cr}$$
and extend $L$\ `log-linearly.' By the same methods as in the proof of
part~(i), we see that $L$\ is convex, that $G\prec L$, and that $L$\ is almost
concave. Finally, we set $H(t) = L(t) \lm L(t)$\ to obtain the result.
\endproof

\proclaim Lemma 5.6.2. Let $G$\ be a \phifunction. Then the following are
equivalent. \itemi $G$\ is almost linear.
\itemii There are \Nfunction s $H$\ and $K$\ satisfying \conditionJ\ such
that $G = H\circ K^{-1}$.
\itemiii There are \Nfunction s $H$\ and $K$\ satisfying \conditionJ\ such
that $G = H^{-1}\circ K$.
\itemiv There are \Nfunction s $H$\ and $K$\ satisfying \conditionJ\ and a
number $c<\infty$\ such that $\invc G(t)/t \le H(t)/K(t) \le c\, G(t)/t $\
for all $t>0$.
\itemv There are strictly monotone weight functions $w_0$\ and
$w_1$ satisfying \conditionL\ and a number $c<\infty$\ such that $ \invc G(t)/t
\le w_0(t) w_1(t) \le c\, G(t)/t $\ for all $t\ge 0$.
\moreproclaim

Before proving this result, we will require a couple of technical lemmas.

\proclaim Lemma 5.6.3. If $G_1$\ and $G_2$\ are equivalent \phifunction s, 
and if
one of $G_1$\ or $G_2$\ satisfies the \Deltacond, then there is a number
$c<\infty$\ such that $\invc G_1(t) \le G_2(t) \le c\, G_1(t)$\ for all $t\ge
0$.

\proclaim Lemma 5.6.4. Suppose that $F:[0,\infty)\to[0,\infty)$\ is a function
such that for some numbers $c_1>1$\ and $c_2>1$\ we have that $F(c_1 t) \ge c_2
F(t)$\ for all $t\ge 0$. Then there is a number $c<\infty$\ and a dilatory
\phifunction\ $G$\ such that $G(\invc t) \le F(t) \le G(c t)$\ for all $t\ge0$.

\Proof: We have that $F(c_1^n t) \ge c_2^n F(t)$\ for all $n\in\N$\ and
$t\ge0$. Then it is clear that
$$ G(t) = \sup_{s\ge1} s^{-\log c_2/\log c_1} F(st) $$
satisfies the conclusion of the lemma.
\endproof

\Proofof Lemma 5.6.2: The implications (ii)\implies(i) and (iii)\implies(i)
follow from Lemmas~5.4.5 and~5.4.6. The implications (iv)\implies(i) and
(v)\implies(i) are obvious.

To show that (i)\implies(ii), we note that, since $G$\ is almost
convex, by Lemma~5.6.1(i), there is an \Nfunction\ $K_0$\ satisfying
\conditionJ\ such that $G\succ K_0^{-1}$. If we let $K(t) = K_0(t\lm t)$, then we
see that $H = G\circ K$\ is an \Nfunction. Since $G$\ is almost concave, it
follows by Lemma~5.4.6 that $H$\ satisfies \conditionJ.
The implication (i)\implies(iii) is similar, using Lemma~5.6.1(ii).

To show (i)\implies(iv), we note that since $G$\ is almost concave, by
Lemma~5.6.1(ii), there is an \Nfunction\ $H_0$\ satisfying 
\conditionJ\ such that
$G\prec H_0$. Then, from Lemmas~5.4.2 and~5.6.4, it follows that $tH_0(t)/G(t)$\
is equivalent to a convex function $K_0$. Since $G$\ is almost convex, we have
that $K_0$\ is almost concave. Now we let $H(t) = H_0(t)\lm t$\ and $K(t) =
K_0(t)\lm t$, and the result follows by Lemma~5.6.3.

To show that (iv)\implies(v), by Lemma~5.4.2, we may assume that
$H$\ and $K$\ are convex. Thus, if we let
$w_0(t) = (\lm t) H(t)/t$\ and $w_1(t) = t/(\lm t )K(t)$, then $\tilde w_0$\ and
$\widetilde{1/w_1}$\ are both almost constant \phifunction s. Then it follows
from Lemma~5.5.1 that $w_0$\ and $w_1$\ satisfy \conditionL.
\endproof

\beginsection 5.7.\ The Proof of the Results in Section 4

Now we are ready to piece together all the lemmas we have just proved.

\Proofof Theorem 4.2:  First we will show that (ii)\implies(i). By Lemma~5.2.1,
we know that $G_2\circ G_1^{-1}$\ is dilatory, and hence if $G_1$\ is dilatory,
then so is $G_2$. Therefore, we may assume that $G_2$\ is dilatory. 

By Lemma~5.3.1, there is a number $c_1<\infty$ such that $\normo
f_{1,G_1\circ G_2^{-1}} \le c_1\, \normo f_1$\ for all measurable $f$. Since
$G_2$\ is dilatory, the result follows by follows by Lemmas~5.1.3 and~5.1.2.

Now we show that (i)\implies(ii). By Lemma~5.2.5, we know that
$G_2$\ satisfies the \Deltacond. Therefore, by Lemmas~5.1.2 and~5.1.3, 
it follows that there is a number $c_1<\infty$\ such
that $\normo f_{1,G_1\circ G_2^{-1}} \le c_1\, \normo f_1$\ for all measurable
$f$. Now the result follows by Lemma~5.3.1.

The implication (ii)\iff(iv) follows from Lemma~5.4.7(i), and (ii)\iff(iii)
follows from Lemma~5.6.1.
\endproof

\Proofof Theorem 4.3:  First we show that (ii)\implies(i). If $G_1\circ
G_2^{-1}$\ is almost linear, then by Lemma~5.2.1, it follows that $G_1\circ
G_2^{-1}$\ satisfies the \Deltacond. Therefore, by Lemma~5.4.3, we see
that $G_2\circ G_1^{-1}$\ is almost convex. Therefore, by Theorem~4.2, we have
that $L_{F_1,G_1}$\ and $L_{F_1,G_2}$\ are equivalent. Clearly, we have
that both $G_1$\ and $G_2$\ are dilatory, and by
Proposition~4.1, $L_{F_1,G_2}$\ and $L_{F_2,G_2}$\ are equivalent. The result
follows.

Next we show that (i)\implies(ii). First notice that, since $\tilde
F_1^{-1}(t) = \normo{\chi_{[0,t]}}_{F_1,G_1}$\ and  $\tilde
F_2^{-1}(t) = \normo{\chi_{[0,t]}}_{F_2,G_2}$, we have that $F_1 \asymp F_2$.
Now let us suppose without loss of generality that $G_1$\ is dilatory. Then by
Proposition~4.1, we have that $L_{F_1,G_1}$\ and $L_{F_2,G_1}$\ are equivalent,
and hence $L_{F_2,G_1}$\ and $L_{F_2,G_2}$\ are equivalent. Now, by
Lemmas~5.2.3 and~5.2.5, both $G_1$\ and $G_2$\ are dilatory, and both $G_1$\ and
$G_2$\ satisfy the \Deltacond. Therefore, $G_2\circ G_1^{-1}$\ and
$G_1\circ G_2^{-1}$\ satisfy the
\Deltacond. By Theorem~4.2, both $G_1\circ G_2^{-1}$\ and $G_2\circ
G_1^{-1}$\ are almost convex. Now the result follows by Lemmas~5.4.2 and~5.4.5.

The implications (ii)\iff(v) follow from Lemma~5.4.8,
the implication (v)\implies(viii) follows from Lemmas~5.6.4 and~5.6.3, and the
implication (viii)\implies(v) is obvious.
Finally, (viii)\iff(iii)\iff(iv)\iff(vi)\iff(vii) all follow from Lemma~5.6.2.
\endproof

\Proofof Theorem 4.4:
The implication (ii)\implies(i) is obvious, so we show that (i)\implies(ii). As
in the proof of Theorem~4.2, we may suppose that $G_2$\ satisfies the
\Deltacond. By Lemma~5.1.2, we may assume without loss of generality that $F(t) =
t$. Now the result follows by Lemmas~5.1.3 and~5.3.1, in the same manner as in
the proof of Theorem~4.2.
\endproof

\Proofof Theorem 4.5:
The implications (iii)\implies(i) and (iii)\implies(ii) are obvious. We
show that (i)\implies(ii). First notice that $\chi_{[0,t)} \in \T_{F_1}$, and
so as $\tilde F_1^{-1}(t) = \normo{\chi_{[0,t)}}_{F_1,G_1}$\ and  $\tilde
F_2^{-1}(t) = \normo{\chi_{[0,t)}}_{F_2,G_2}$, we have that $F_1$\ and $F_2$\
are equivalent. Then it is clear that there is a number $c<\infty$\ such that
$f^*\in\T_{F_1}$\ if and only if there is a function $g^*\in\T_{F_2}$\ such that
$\invc f^* \le g^* \le c\,f^*$. Similarly, (ii)\implies(i).

Now we show that (i) and (ii)\implies(iii). Let us suppose, without loss of
generality, that $G_1$\ is dilatory. Then by
Proposition~4.1, we have that $L_{F_1,G_1}$\ and $L_{F_2,G_1}$\ are equivalent,
and hence $L_{F_2,G_1}$\ and $L_{F_2,G_2}$\ are equivalent. Now, by Lemma~5.2.3,
it follows that if one of $G_1$\ or $G_2$\ is dilatory, then both are. Then the
result follows by Theorem~4.4. 
\endproof

\Proofof Theorem 4.6: We show that (i)\iff(iii). By Proposition~4.1, we know
that if $F_1$\ and $F_2$\
are equivalent, then $L_{F_1,\infty}$\ and $L_{F_2,\infty}$\ are
equivalent. Also, if $L_{F_1,G}$\ and $L_{F_2,\infty}$\ are
equivalent, then since $\tilde F_1^{-1}(t) =
\normo{\chi_{[0,t]}}_{F_1,G}$\ and  $\tilde F_2^{-1}(t) =
\normo{\chi_{[0,t]}}_{F_2,\infty}$, we have that $F_1$\ and $F_2$\
are equivalent. Thus,
without loss of generality, we may assume that $F_1 = F_2 = F$.

Now we note that we always have that $\normo f_{F,\infty} \le \normo
f_{F,G}$. This follows because for all $x\ge0$, we have that $f^* \ge f^*(x)
\chi_{[0,x]}$, and hence 
$$ \normo f_{F,G} \ge \normo{f^*(x)\chi_{[0,x]}}_{F,G} \ge 
   f^*(x) \tilde F^{-1}(x) .$$

Now we show that $L_{F,G}$\ is equivalent to $L_{F,\infty}$\ if and only if
$\normo{1/\tilde F^{-1}}_{F,G} < \infty$. That the first statement implies the
second is obvious, because $\normo{1/\tilde F^{-1}}_{F,\infty} = 1$. To show the
second statement from the first, note that if $\normo f_{F,\infty} \le 1$, then
$f^*(x) \le 1/\tilde F^{-1}$.

But $\normo{1/\tilde F^{-1}}_{F,G} = \normo{1/\tilde G^{-1}}_G$, and the result
follows. The other implications follow by Lemma~5.5.2.
\endproof

\beginsection 6.\enspace Is Every R.I.\ Space
Equivalent to an Orlicz--Lorentz Space?

We can answer the question in the negative easily, as follows. It is well known
that $L_{1,\infty}$\ is not separable. Then it is not hard to see that
$L_{1,\infty}^0$, the closure of the simple functions in $L_{1,\infty}$, is not
an Orlicz--Lorentz space.

However, the reader may consider this cheating. So to avoid all this `infinite
dimensional nonsense,' we might ask the following question. Is there a
rearrangement invariant space $X$\ such that for all Orlicz--Lorentz spaces
$L_{F,G}$, the norms $\normdot_X$\ and $\normdot_{F,G}$\ are inequivalent on
$X\cap L_{F,G}$? We answer this question in the positive by the following
example.

\proclaim Theorem 6.1. There is a rearrangement invariant Banach space $X$, where
the measure space is $[0,\infty)$\ with Lebesgue measure, such that for every
Orlicz--Lorentz space $L_{F,G}$, the norms $\normdot_X$\ and $\normdot_{F,G}$\
are inequivalent on the vector space of simple functions.

\Proof: We define the
following norm for measurable functions $f$:
$$ \normo f_X = \sup\{\, \normo{f^* g}_1/\normo g_2 : g \in \T_{t^2} \} .$$
We let $X$\ be the vector space of all measurable functions $f$\ such that
$\normo f_X < \infty$, modulo functions that are zero almost everywhere. Then
it is an easy matter to see that $X$\ is a rearrangement invariant space such
that $\normo
g_X = \normo g_2$\ for all $g\in\T_{t^2}$.
Thus, if for some \phifunction s $F$\ and $G$\ we have that
$$ c_1^{-1} \normo f_X \le \normo f_{F,G} \le c_1 \normo f_X $$
for all simple functions $f$, then by Theorem~4.5, we see that there is a
constant $c_2<\infty$\ such that
$$ c_2^{-1} \normo f_X \le \normo f_2 \le c_2 \normo f_X $$
for all simple functions $f$. We will show that this cannot happen.

Define $\T'$\ to be the set of functions $h:[0,\infty)\to[0,\infty)$\
such that for some integers $k_1 < \ldots < k_n$, and setting $k_0=-\infty$, we
have that 
$$ h(x) = \cases{
   2^{-k_i} & if $4^{k_{i-1}}\le x < 4^{k_i}$\ and $1\le i\le n$ \cr
   0        & otherwise. \cr }$$
Then it is an easy matter to see that if $g\in \T_{t^2}$, then there is a
function $h\in \T'$\ such that $h(4x)/2 \le g(x) \le 2h(x/4) $. Therefore, we
see that
$$ \ts{1\over 4} \normo f_X
   \le \sup\{\, \normo{f^* h}_1 / \normo h_2 : h\in\T' \,\}
   \le 4 \normo f_X .$$

Now, for each $N\in\N$, let $f_N$\ be the simple function
$$ f_N(x) = 
   \cases{ 1               & if $0\le t < 4$ \cr
           k^{-1/2} 2^{-k} & if $4^k \le t < 4^{k+1}$\ and $1\le k \le N$\cr
           0               & otherwise. \cr} $$
Then it is easy to see that $\normo{f_N}_2 \to \infty$\ as $N\to \infty$.
However, a simple, but laborious, calculation shows that there is a number
$c<\infty$\ such that $\normo{f_N^* h}_1/\normo h_2 \le c$\ for all $h\in\T'$,
and hence $\normo{f_N}_X \le 4c$. 
\endproof

\beginsection 7.\enspace The Definition of Torchinsky and Raynaud

The definition of the Orlicz--Lorentz spaces presented here is not the only
possible definition. In fact, given any weight function
$w$\ and any \phifunction s $H$\ and $G$, one can form
the functional $$ \normo f_{w,H,G} = \normo{w \cdot (f^* \circ H)}_G .$$
We have investigated the case when $w(x) = 1$. However, Torchinsky \ref[T]
gave the following definition for the Orlicz--Lorentz functional. If $F$\ and
$G$\ are \phifunction s, then we define 
$$ \eqalignno{
   \normo f_{F,G}^T
   &= \normo{\tilde F^{-1}(e^x) f^*(e^x)}_G \cr
   &= \inf \left\{\, c :
      \int_0^\infty G\bigl(\tilde F^{-1}(x) f^*(x)/c\bigr)
      \, {dx\over x} \le 1 \right\} ,\cr}$$
and call the corresponding space $L_{F,G}^T$\ (my notation).

These spaces were investigated by Raynaud \ref[R]. He showed that if $F$\ is
dilatory and satisfies the \Deltacond, and if $G$\ is dilatory, then
$\normo{\chi_A}_{F,G}^T \approx \tilde F^{-1}\bigl(\mu(A)\bigr)$\ for all
measurable $A$. Thus, it follows that
$L_{F,G}^T$\ and $L_{F,G}$\ are equivalent if $G(t) = t^p$.

The comparison results for these spaces are much more straightforward. Raynaud
\ref[R] showed that if $F_1$\ and $F_2$\ are
dilatory and satisfy the \Deltacond, and if $G_1$\  and $G_2$\ are dilatory,
then $L_{F_1,G_1}^T$\ and $L_{F_2,G_2}^T$\ are equivalent if
$F_1$\ and $F_2$\ are equivalent, and the {\it sequence\/} spaces $l_{G_1}$\ and
$l_{G_2}$\ are equivalent. The converse result is also easy to show.

We also comment that the Boyd indices of these spaces are much easier
to compute. This will be dealt with more fully in \ref[Mo2].

Also, unlike the Orlicz--Lorentz spaces we have used here, we do not always have
that $L_{F,F}^T$\ is equivalent to the Orlicz space $L_F$. For example, if $F(t)
= t \lm t$, then $L_F$\ is equivalent to $L_{F,1}^T$\ by Theorem~4.3, and since
$l_F$\ and $l_1$\ are not equivalent, this is not equivalent to $L_{F,F}^T$.

We finally add that we may define the spaces $L_{F,X}$, where $X$\ is a
rearrangement invariant quasi-Banach space on $\R$\ satisfying certain mild
restrictions. Corresponding to the definition used in this paper, we may define
$$ \normo f_{F,X} = \normo{f^*\circ \tilde F \circ \phi_X}_X ,$$
where $\phi_X(t) = \normo{\chi_{[0,t]}}_X$\ is the fundamental function of $X$.
Corresponding to the definition used by Torchinsky and Raynaud, we may define
$$ \normo f_{F,X}^T = \normo{\tilde F^{-1}(e^x) f^*(e^x)}_X .$$

\beginsection Acknowledgements

This paper is an extension of work that I presented in my Ph.D.\ thesis
\ref[Mo1]. I would like to express my thanks to D.J.H.~Garling, my Ph.D.\
advisor, as well as the Science and Engineering Research Council who financed my
studies at that time.
I would also like to express much gratitude to A.~Kami\'nska, W.~Koslowski and
N.J.~Kalton for their keen interest and useful conversations. I also
extend my thanks to the referee for his diligent reading of the final 
manuscript (although all errors are, of course, my responsibility). 
Finally, I would
like to mention T.N.~Prokraz for his great influence on this work.

\beginsection References

\references{
B--R & C.~Bennett and K.~Rudnick,\rm\ On Lorentz--Zygmund spaces,\sl\
Dissert.\
Math.\ {\bf 175} (1980), 1--72.\cr
B--S & C.~Bennett and R.~Sharpley,\sl\ Interpolation of Operators,\rm\
Academic Press 1988.\cr
H & R.A.~Hunt,\rm\ On $L(p,q)$\ spaces,\sl\ L'Enseignement Math.\ (2)
{\bf 12} (1966), 249--275.\cr
Ka1 & A.~Kami\'nska, \rm\ Some remarks on Orlicz--Lorentz spaces,\sl\ Math.\
Nachr., to appear.\cr
Ka2 & A.~Kami\'nska, \rm\ Extreme points in Orlicz--Lorentz spaces,\sl\
Arch.\ Math., to appear.\cr
Ka3 & A.~Kami\'nska, \rm\ Uniform convexity of generalized Lorentz spaces,\sl\
Arch.\ Math., to appear.\cr
K--R & M.A.~Krasnosel'ski\u\i\ and Ya.B.~Ruticki\u\i,\sl\
Convex Functions and
Orlicz Spaces,\rm\ P.~Noodhoof Ltd., 1961.\cr
Lo1 & G.G.~Lorentz,\rm\ Some new function spaces,\sl\ Ann.\ Math.\ {\bf 51}
(1950), 37--55.\cr
Lo2 & G.G.~Lorentz,\rm\ On the theory of spaces $\Lambda$,\sl\ Pac.\ J.\ Math.\
{\bf 1} (1951), 411--429.\cr
Lo3 & G.G.~Lorentz,\rm\ Relations between function spaces,\sl\ Proc.\ A.M.S.\
{\bf 12} (1961), 127--132.\cr
Lu & W.A.J.~Luxemburg,\sl\ Banach Function Spaces,\rm\ Thesis, Delft Technical
Univ., 1955.\cr
Ma & L.~Maligranda,\rm\ Indices and interpolation,\sl\ Dissert.\ Math.\ {\bf 234}
(1984), 1--49.\cr
My & M.~Masty\l o,\rm\ Interpolation of linear operators in
Calderon--Lozanovskii spaces,\sl\ Comment.\ Math.\ {\bf 26,2} (1986),
247--256.\cr
Mo1 & S.J.~Montgomery-Smith,\sl\ The Cotype of Operators from
$C(K)$,\rm\
Ph.D.\ thesis, Cambridge, August 1988.\cr
Mo2 & S.J.~Montgomery-Smith,\rm\ Boyd Indices of Orlicz--Lorentz spaces,\sl\
in preparation.\cr
O & W.~Orlicz,\rm\ \"Uber eine gewisse Klasse von R\"aumen vom Typus B, \sl\
Bull.\ Intern.\ Acad.\ Pol.\ {\bf 8} (1932), 207--220.\cr
R & Y.~Raynaud,\rm\ On Lorentz--Sharpley spaces,\sl\ Proceedings of the
Workshop ``Interpolation Spaces and Related Topics'', Haifa, June 1990.\cr
S & R.~Sharpley,\rm\ Spaces $\Lambda_\alpha(X)$\ and Interpolation,\sl\ J.\
Funct.\ Anal.\ {\bf 11} (1972), 479--513.\cr
T & A.~Torchinsky,\rm\ Interpolation of operators and Orlicz classes,\sl\
Studia Math.\ {\bf 59} (1976), 177--207.\cr
}

\bye